\documentclass[twoside, 12pt]{article}
\usepackage[cp1251]{inputenc}
\usepackage[russian, english]{babel}
\usepackage{epsfig}
\usepackage{amssymb,amsmath,amsfonts,soul,amsthm,enumerate}
\usepackage{amsfonts}
\usepackage{amsmath}
\usepackage{graphicx}
\usepackage{amsthm}
\usepackage[T2A]{fontenc}
\usepackage{mathrsfs}
\usepackage[english]{babel}
\newtheorem{theorem}{Theorem}[section]

\newtheorem{definition}{Definition}
\newtheorem{remark}{Remark}

\numberwithin{equation}{section}
 \def\@evenhead{\vbox{\hbox to \textwidth{\thepage\hfil\sl\leftmark\strut}\hrule}}
 \def\@oddhead{\vbox{\hbox to \textwidth{\rightmark\hfill\thepage\strut}\hrule}}

\begin{document}
 \sloppy

\centerline{\bf SHARP RESTRICTIONS OF ANALYTIC FUNCTION SPACES}
\centerline{\bf OF SEVERAL VARIABLES} 

\vskip 0.3cm

\centerline{\bf R.F.~Shamoyan, N.M.~Makhina} 

\markboth{\hfill{\footnotesize\rm   R.F.~Shamoyan, N.M.~Makhina  }\hfill}
{\hfill{\footnotesize\sl  Sharp restrictions of analytic function spaces of several variables}\hfill}
\vskip 0.3cm

\vskip 0.7 cm


\vskip 0.2cm


\vskip 0.2cm

\noindent {\bf Abstract.} In this expository paper we collect many recent advances in analytic function spaces of several complex variables related with trace problem in tubular domains over symmetric cones and bounded strongly pseudoconvex domains with smooth boundary. We consider various function space of analytic functions of several variables in various domains in $\mathbb C^n$ and provide or complete descriptions of traces or estimates of traces of various analytic function spaces in various domains  obtained in recent years, by various authors.
The problem to find sharp estimates of traces of Hardy analytic function spaces in the unit polydisk first was posed by W.~Rudin in 1969. Since then many papers appeared in literature. We collect in this expository paper not only already many  known results on traces of various analytic  function spaces in product domains  but also discuss various new interesting results related with this problem. Related to trace problem various results were provided previously by G.~Henkin, E.~Amar, H.~Alexander and various other authors. Finnaly, note that our trace theorems are closely related with the Bergman type projections acting between function spaces with different dimensions. This expository paper contains mainly new results concerning traces in tubular and bounded strongly pseudoconvex domains, proofs of these theorems are based in particular also on various properties of Bergman type projection, in this expository paper we will also shortly  discuss  some new results obtained by first author on Bergman type projections in these complicated domains in $\mathbb C^n.$ This is the second part of our notes related with trace problem. In the first part we provided a large list of recent sharp results on traces in the polydisk and polyball and was published in \cite{ShMakh2025}. 
In the last section of this expository paper we collect many interesting remarks and put various open interesting questions for interested readers and also add various interesting short comments related with this problem of traces in various simlpe and complicated domains in $\mathbb C^n$ in various analytic function spaces of several variables. Some new results concerning traces in harmonic function spaces of several variables will be also added.

\section{\large Introduction}

This paper is the second part of our notes related with trace problem in analytic function spaces of several complex variables, namely in tube, pseudoconvex and related domains. The first part provided large list of sharp recent results in the polydisk and unit polyball and was published recently in \cite{ShMakh2025}.

Let $G$ be any bounded or unbounded domain in $\mathbb{C}^n$. Let $f$ be an analytic function in $G$ from a certain fixed quazinormed analytic function space $X(G)$ in $G$. This expository paper concerns the following problem to find precise estimates of $|f(z,...,z)|$ set where $z$ is an arbitrary point from $G$ in ${\mathbb C}^n$.

Obviously, very similarly such a problem may be posed in tubular domains over symmetric cones and bounded strongly pseudoconvex domains with smooth boundary and in products of such type domains in ${\mathbb C}^n$.

This problem was posed for the first time in W. Rudins book \cite{Rud1969} in 1969 for traces of Hardy spaces in bidisk and since then many papers appeared in this direction. We mention for example old papers \cite{Dur1975}-\cite{ShF1990} and various references there.

We started \cite{ShMakh2025} with Trace results first in the simplest domain namely the unit polydisk and the unit polyball. In this paper we provide new results. 

One of the goals of this paper is to consider a new trace map in the trace problem related with embedding theorems for the case of the bounded strongly pseudoconvex (or tubular) domains $D$ with a smooth boundary; it is a map $T_rf(z) =f(z, . . . , z),$ $z \in D,$ for $f \in X \subset H(D^m)$ for a certain quasinormed analytic space $X$ on $D^m.$ 

We give precise definitions. 

\begin{definition}\label{def1} {\rm Let D be bounded or unbounded domain in $\mathbb C^n$. Let $H(D^m)$ be the space of all analytic functions in $D^m$, $D^m = D\times. . .\times D,$ $m \in \mathbb{N}$. We say $Trace (X) = Y,$ if $X$ is a quasinormed subspace of $H(D^m)$ and $Y$ is a certain fixed quasinormed subspace of $H(D),$ if for each $F,$ $F \in X,$ $F(z, . . . , z) = f(z),$ $f \in Y,$ and the reverse is also true: for each $g,$ $g \in Y,$ there is a function $F,$ $F \in X,$ such that $F(z, . . . , z) = g(z),$ for all $z \in D.$}
\end{definition}

This definition of traces of analytic function  spaces clearly may be easily extended to function spaces consisting of harmonic and pluriharmonic functions in product domains of several variables.

We note that we replace $Trace (X) = Trace X$ below by simpler $DX$ sometimes. 

Various new interesting estimates, and sharp results were obtained in this direction in recent years. See, for example, recent papers \cite{Sh2007V}-\cite{ShM2022R} and also many references there also.

The goal of this paper to collect all these new interesting results in this new expository research paper. Note that various new open interesting problems related with this problem will also be provided in this paper.

In this expository paper we will collect many old and new results concerning this problem for many X analytic function spaces in various  product domains in ${\mathbb C}^n$. First below we provide definitions needed for formulations of such type results, then we formulate our results (old and new) on Traces. 

\begin{remark} {\rm We also note that very similar problems, Trace theorems or various estimates of traces of functions were given, proved and considered in many  complicated functional spaces in ${\mathbb R}^n$ by many authors. We refer  the reader, for example, to \cite{Bes1975}, \cite{Maz1985}, \cite{Trieb1992} and many references there.}
\end{remark}

Trace theorems have also interesting extensions and also nice applications; we, for example, refer the reader to the papers of E. Amar, C. Menini \cite{Am2002} and D. Clark \cite{Cl1988} where such type applications and such type extentions were provided.

Finally note that  our trace theorems are closely related with the Bergman type projections acting between function spaces with different dimensions.

For the convenience of the reader we organize this expository  paper in the following simple manner. We divide all results concerning traces and related results  by domains. Namely various Theorems on traces concerning different domains will be collected in different sections.

The special attention in this paper will be given to recent research papers of the first author related with the trace problem and the plan of the paper is the following: we provide first various definitions of various interesting analytic function spaces in various domains in $\mathbb C^n$ then provide some precise information on traces of such type function spaces in $\mathbb C$. Note also that the trace operators may be also considered in relation with various embedding theorems and such type sharp results will also be given in this expository paper very shortly in tubular or pseudoconvex domains.

In the last section of this expository paper, we collect many interesting remarks and put various open interesting questions for interested readers and also add various interesting short comments related to this problem of traces in various simple and complicated domains in $\mathbb C^n$ in various analytic function spaces of several variables.

Some new results concerning traces in harmonic function spaces of several variables will also be added.

The interested reader may easily note that in this paper all trace theorems and analytic spaces in tubular domains have complete analogues in pseudoconvex domains and formulations of our results concerning traces in these very different domains in $\mathbb{C}^n$ also have similarities.

Note in \cite{Sh2007V}-\cite{ShM2010S} the first author started his intensive investigations on traces in various analytic functions in simple domains as bidisk, polydisk, polyball. And later these issues were considered by him in more complicated domains in $\mathbb{C}^n$ as tube and bounded strongly pseudoconvex domains and analytic function spaces of mixed norm type and Bergman type on them.

\section{\large Traces in various spaces of analytic functions in strongly pseudoconvex domains}

In  this section we collect several recent results on traces in the polyball and pseudoconvex domains from papers of the first author \cite{ShM2015J}, \cite{ShM2015DAN}, \cite{ShM2017C}. We define below various analytic function spaces in the context of bounded pseudoconvex domains and provide descriptions of their traces or estimates of sizes of their traces. 

First the first author provided many results (sharp results) in the unit ball on traces with O.~Mihic, then later they were extended to similar type analytic function spaces in bounded strongly psrudoconvex domains with smooth boundary in $\mathbb{C}^n$ in his research papers with O.~Mihic and S.~Kurilenko.

We refer for the differential $\mathcal{D}^{\alpha}$ operator which appears below in definition of our spaces to paper of J.~Fabrega and J.~Ortega \cite{OrtFab1994} which is mentioned there  by us. 

We refer the interested reader for some basic information on Martinelly--Bochner integrals and $U(z,w)$ Martinelly--Bochner  kernels which appear below to research papers which are mentioned in those theorems below.

Let  $D\subset {\mathbb{C}}^n$ be a domain, that is, an open connected subset. One says that $D$ is pseudoconvex (or Hartogs pseudoconvex) if there exists a continuous plurisubharmonic function $\varphi$ on $G$ such that the set  $\{ z \in D \mid \varphi(z) < x \}$ is a relatively compact subset of $D$ for all real numbers $x.$ In other words, a domain is pseudoconvex if $D$ has a continuous plurisubharmonic exhaustion function. Every (geometrically) convex set is pseudoconvex. However, there are pseudoconvex domains which are not geometrically convex.

When $D$ has a $C^2$ (twice continuously differentiable) boundary, this notion is the same as Levi pseudoconvexity, which is easier to work with. More specifically, with a $C^2$ boundary, it can be shown that  $D$ has a defining function, i.e., that there exists  $\rho: \mathbb{C}^n \to  \mathbb{R}$ which is $C^2$ so that $D=\{\rho <0 \},$ and $\partial D =\{\rho =0\}.$  Now,  $D$ is pseudoconvex iff for every  $p \in \partial G$ and $w$ in the complex tangent space at $p,$ that is,  

\[\nabla \rho(p) w = \sum_{i=1}^n \frac{\partial \rho (p)}{ \partial z_j }w_j =0,\] we have
\[\sum_{i,j=1}^n \frac{\partial^2 \rho(p)}{\partial z_i \partial \bar{z}_j } w_i \bar{w}_j \geq 0.\]

The definition above is analogous to definitions of convexity in real analysis.

Let $D$ is the pseudoconvex domain in $\mathbb C$,  $\tilde{\mathcal{D}}^{\alpha } $ - differential operator in $D,$ $ 0<p,q<+\infty,$ Hardy-Sobolev and Bergman-Besov spaces of analytic functions
\[H_{\alpha ,\beta }^{p} (D)=\left\{f\in H(D):\mathop{\sup }\limits_{\varepsilon >0} \left(\int _{\partial D_{\varepsilon } }\left|\tilde{\mathcal D}^{\alpha } f(\zeta )\right|^{p} d\sigma _{\zeta }  \right)^{\frac{1}{p} } \varepsilon ^{\beta } <+\infty \right\},\beta \ge 0,\alpha >0.\] 

Denote (see \cite{OrtFab1994})\\
$A_{\delta ,k}^{p,q} (D)=\Bigg\{f\in H(D):$
\[\sum _{\left|\alpha \right|\le k}\int _{0}^{r_{0} }\left(\int _{\partial D_{r} }\left|\tilde{\mathcal D}^{\alpha } f(\zeta )\right|^{p} d\sigma _r  \right)^{\frac{q}{p} }   r^{\delta \frac{q}{p} -1} dr<+\infty \Bigg\},\delta >0,\alpha >0,\]
where  $D_{r}=\left\{z \in \mathbb{C}^{n}: \rho(z)<(-r)\right\},$ $\partial D_{r}$ is a boundary, $d \sigma_{r}$ is a normalized surface measure on $\partial D_{r}$ and by $dr$ we denote a normalized volume element on $(0, r),$ $0<p<\infty,$ $0<q \leq \infty,$ $\delta>0,$ $k=0,1,2 \ldots$ and
\[\|f\|_{p, \infty, \delta, k}=\sup \left\{\left(\sum_{|\alpha| \leq k \mid}\left(r^{\delta}\right) \int_{\partial D_{r}}\left|\mathcal{D}^{\alpha} f\right|^{p} d \sigma_{r}\right)^{\frac{1}{p}}: 0<r<r_{0}\right\},\]
(for $p, q<1$ it is quazinorm) (see \cite{OrtFab1994}). 
It can be easily shown these spaces are the Banach spaces for $\min (p, q) \geq 1$. Further, for $p=q$ we have
\[\|f\|_{p, \delta, k}=\left(\sum_{|\alpha| \leq k} \int_{D}\left|\mathcal{D}^{\alpha} f(\zeta)\right|^{p}(-\rho(\zeta))^{\delta-1} d v(\zeta)\right)^{\frac{1}{p}} ; \delta>0, k \geq 0,\]
where $d v(\zeta)$ (or sometimes $d V(\zeta)$ ) is a normalized Lebegues measure on $D$ (for $k=0$ we will use another, more convenient, notation for $A_{\alpha}^{p}$ spaces). As we see for $p=q;$ $ k=0$, we get the usual Bergman spaces $A_{\delta}^{p}(D)=A_{\alpha}^{p}(D)$.

We define same spaces in product domains in a natural way, $D^m = D\times \dots \times D.$ 

In theorem \ref{th2.1}-\ref{th2.2} provide some analogues of results on traces in terms of Martinelly--Bochner integrals and kernels. 

We formulate theorem \ref{th2.1} in context of unit ball $B$ in $\mathbb {C}^n$ but similar results may be valid in bounded strongly pseudoconvex domains with smooth boundary.

\begin{theorem}\label{th2.1} {\rm (see \cite{ShK2015C}).}
\begin{enumerate}
\item  Let $p\le 1,\beta \ge \alpha \ge 0,$ then $A_{nm+(\beta -\alpha )pm-n}^{p,p} (B)\subset Trace(H_{\alpha ,\beta }^{p} (B^{m} ));$
\item  Let $p\le 1,$ $\displaystyle \alpha _{j} \ge 0,t\ge \frac{\sum _{j=1}^{m}\alpha _{j}  }{m} ,$ $\beta _{j} =\alpha _{j} +1,j=1,...,m.$ 

Then $A_{mt+mn-n-\sum _{j=1}^{m}\alpha _{j}  }^{p,p} (B)\subset Trace(A_{t,\vec{\beta }}^{p,p} (B^{m} )).$ 
\end{enumerate}
\end{theorem}

\begin{theorem} \label{th2.2} {\rm (see \cite{ShK2015C}).}
\begin{enumerate}
\item  Let $p\le 1,\alpha \ge \beta ,$ then $Trace(H_{\alpha ,\beta }^{p} (B^{m} ))$ contains a f function for which 
 \[f(z)=\int _{\partial D}f(\zeta )U(\zeta ,z);z\in D. \eqno(IR)\]
\item  Let $p\le 1,$ $\displaystyle \alpha _{j} \ge 0,$ $t\ge \frac{\sum \limits _{j=1}^{m}\alpha _{j}  }{m} ,$ $\beta _{j} =\alpha _{j} +1,j=1,...,m.$ 

Then $Trace(A_{t,\vec{\beta }}^{p} (B^{m} ))$ contains a $f$ function for which (IR) integral representation is valid.
\end{enumerate}
\end{theorem}

In \cite{ShK2015I} considered analytic Bloch spaces in strongly pseudoconvex domains with smooth boundary, spaces of mixed norm and the so-called new spaces of Hertz type of analytic functions in domains of the same type.

Let $D$ be a bounded strongly pseudoconvex domain of ${\mathbb C}^{n} $ with $C^{\infty } $ boundary. We assume that the
strongly plurisubharmonic function $\rho$  is of class $\mathbb C^{\infty}$ in a neighborhood of $\bar {D}$, that $-1 \le \rho <0,$ $z \in D,$ and $|\partial \rho | \ge c_0 > 0$ for $|\rho | \ge r_0.$ 

Let also $\delta :D\to \mathbb R^{+} $ will this Euclidian distance from the boundary, that is $\delta (z)=d(z,\partial D)$; given $r\in (0,1),z_{0} \in D,$ we shall denote by $B_{D} (z_{0} ,r)$ ball of center $z_0$ and radius $\displaystyle \frac{1}{2} \log \frac{1+r}{1-r} $. 

We shall also need the existence of suitable coverings by Kobayashi balls:

\begin{definition} {\rm Let $D \subset \mathbb{C}^{n}$ be a bounded domain, and $r>0$. An $r$-lattice in $D$ is a sequence $\left\{a_{k}\right\} \subset D$ such that $D=\bigcup_{k} B_{D}\left(a_{k}, r\right)$ and there exists $m>0$ such that any point in $D$ belongs to at most $m$ balls of the form $B_{D}\left(a_{k}, R\right)$, where $R=\frac{1}{2}(1+r)$.}
\end{definition}

Let $\Delta _{k} =\left|B(a_{k} ,r)\right|$, $d\upsilon (z_{1} ,...,z_{m} )=\prod _{j=1}^{m}d\upsilon (z_{j} )=d\upsilon (\vec{z}),z_{j} \in D,j=1,...,m, $ be the normalized Lebeques measure on a product domain. We denote, as usual, by $d\upsilon _{\gamma } =\delta ^{\gamma } (z)d\upsilon (z)$ the weighted Lebegues measure on D and similarly on products of such domains using products of $\delta $ functions in a standard way. We define new analytic Herz-type spaces on products of pseudoconvex domains with a smooth boundary in ${\mathbb C}^{n} $ as follows. 

Let it further be 
\[K_{\alpha ,\beta }^{p,q} (D^m)=\{f\in H(D^{m} ):\] 
\[\int _{D}...\int _{D}\left(\int _{B(\tilde{z}_{1} ,R)}...\int _{B(\tilde{z}_{1} ,R)}\left|f(z_{1} ,...,z_{m} )\right|^{p} \delta ^{\alpha } (z_{1} )...\delta ^{\alpha } (z_{m} )d\upsilon (z_{1} ,...,z_{m} )  \right)^{\frac{q}{p}} \times   \] 
\[\times \prod _{k=1}^{m}\delta ^{\beta } (\tilde{z}_{k} ) d\upsilon (\tilde{z}_{1} ,...,\tilde{z}_{m} )<+\infty ,0<p,q<+\infty ,\alpha >-1,\beta >-1.\}\] 

These are the Banach spaces for $\min (p, q) \geq 1$ and complete metric spaces for the other values of $p$ and $q$. We fix $r$-lattice $\left\{a_{k}\right\} \in D$. For the same values of parameters, we define another analytic Herz type space in pseudoconvex domains with a smooth boundary as follows:

\[K_{\beta}^{p, q}\left(D^{m}\right)=\{f \in H\left(D^{m}\right):\]
\[\sum_{k_{1}} \cdots \sum_{k_{m}}\left(\int_{B(a_{k_1}, r)}\cdots \int_{B(a_{k_m},r)}\left|f\left(z_1, \ldots z_m\right)\right|^p \delta^{\beta}(z_{1}) \ldots
\delta ^{\beta } (z_{m} )d\upsilon (z_{1} ,...,z_{m} )  \right)^{\frac{q}{p}}<+\infty,\]
\[0<p,q<+\infty, \beta >-1\}.\]

If $p=q$, we get the classic Bergman spaces on the products of the pseudoconvex domains with a smooth boundary. This simple fact is based on the properties of r-lattices. These new classes of analytic functions are the Banach spaces for all values of $p, q$ so that $\min (p, q)>1$ and complete metric spaces for other values of parameters. Looking at quazinorms above, we can replace the outer product integration or product sum by one integral or one sum defining similar spaces with the following quazinorms:

\[
\sum_{k}\Big(\int_{B(a_{k},r)} \ldots \int_{B(a_{k},r)} |f(z_1 , \ldots z_m) |^{ p } \delta^{\beta}(z_{1}) \ldots \delta^{\beta} (z_{m})\times \]
\[\times d\upsilon(z_{1}, \ldots, z_{m})\Big)^{\frac{q}{p}}<\infty, 0<p,q<\infty, \beta>-1, \]
\[\int_{D}\Big(\int_{B_{D}(\tilde{z}, R)} \cdots \int_{B_{D}(\tilde{z}, R)}|f(z_{1}, \ldots, z_{m})|^{p} \delta^{\alpha}(z_{1}) \ldots \delta^{\alpha}(z_{m}) \times \]
\[\times d \upsilon(z_{1}, \ldots, z_{m})\Big)^{\frac{q}{p}} \delta^{\beta}(\tilde{z}) d\upsilon(\tilde{z})<+\infty, p,q \in(0, +\infty), \alpha>-1, \beta>-1.\]

Traces of $K_{t,\beta }^{q,p} (D^{m})$ classes can be described in a following manner.  

\begin{theorem} {\rm (see \cite{ShK2015I}).} Let $0<p<+\infty ,$ $t_{j} >-1,$ $\beta _{j} >-1,$ $j=1,...,m,$ $\alpha >-1,$  $\alpha _{2} =\sum\limits _{j=1}^{m}(\beta _{j} +2(n+1)+t_{j} )-(n+1). $  Then $Trace(K_{t,\beta }^{p,p} (D^{m} ))=A_{\alpha _{2} }^{p} $. 

For the same values of parameters we have $Trace(K_{t,\beta }^{p,p} (D^{m} ))=A_{\alpha _{3} }^{p} $, where $\alpha _{2} =\sum\limits _{j=1}^{m}\beta _{j} +(n+1)(m-1). $ 

Moreover, we have $Trace(K_{\beta }^{q,p} (D^{m} ))\subset A_{\alpha _{1} }^{p} ,$ $Trace(K_{\beta _{t} }^{q,p} (D^{m} ))\subset A_{\alpha }^{p} ,$ where $\alpha _{1} =\sum\limits _{j=1}^{m}(\beta _{j} +n+1)p/q-(n+1), $ $\alpha _{1} =\sum\limits _{j=1}^{m}(\beta _{j} +n+1)p/q+\sum\limits _{j=1}^{m}(t_{j} +n+1) -(n+1), $ for all positive values of $p$ and $q$.
\end{theorem}

\begin{remark} {\rm For $m = 1$ some of these assertions are obvious. For the unit ball case in $\mathbb C^n$, these assertions were proved in \cite{ShM2008R}, \cite{ShM2009A}, \cite{ShM2010J}, \cite{ShM2010S}.}
\end{remark}

Below, we provide another new sharp trace theorem for the Bloch type analytic function spaces in the pseudoconvex domains with a smooth boundary in $\mathbb C^n$. 

Let also 
\[A_{\log } (r_{1} ,...,r_{m} )=\left\{f\in H(D^{m} ):\mathop{\sup }\limits_{z_{j} \in D} \left|f(z_{1} ,...,z_{m} )\right|\prod _{j=1}^{m}\log \left(\frac{c}{\delta (z_{j} )} \right)^{-\frac{1}{r_{j} } } \delta ^{\frac{1}{r_{j} } } (z_{j} )<+\infty  \right\},\]
$\displaystyle \sum \limits_{j=1}^{m}\frac{1}{r_{j} } =1,$ $r_{j} >0,$ $c=diam(D);$

\[A_{\log } (1)=\left\{f\in H(D):\mathop{\sup }\limits_{z\in D} \left|f(z)\right|\delta ^{r} (z)\left(\log \frac{c}{\delta (z)} \right)^{-1} <+\infty \right\};\]

\[A_{\vec{r }}^{\infty } (D^{m} )=\left\{f\in H(D^{m} ):\mathop{\sup }\limits_{z_{j} \in D} \left|f(z_{1} ,...,z_{m} )\right|\prod \limits_{j=1}^{m}\delta ^{r_{j} } (z_{j} )<+\infty \right\},\]
$r_{j} >0,$ $j=1,...,m;$
\[A_{r}^{\infty } (D)=\left\{f\in H(D):\mathop{\sup }\limits_{z\in D} \left|f(z)\right|\delta ^{r} (z)<+\infty \right\},\]
$r>0.$

In the following theorem we provide a sharp result on traces of $A_{\vec{r}}^{\infty } (D^{m} )$ and $A_{\log } (r_{1} ,...,r_{m} )$ analytic spaces.

\begin{theorem} {\rm  (see \cite{ShK2015I}).} 
\begin{enumerate}
\item  Suppose $r_{j} >0,\sum\limits _{j=1}^{m}r_{j} =r,j=1,...,m. $ Then $Trace(A_{\vec{r}}^{\infty } (D^{m} ))=A_{r}^{\infty } (D).$
\item  Suppose $r_{j} >0,\sum\limits _{j=1}^{m}r_{j} =1,j=1,...,m. $ Then $Trace(A_{\log } (r_{1} ,...,r_{m} ))=A_{\log } (1).$
\end{enumerate}
\end{theorem}

Let us consider further the traces of Bergman spaces of with the mixed norm on $D^m$:\\
$A_{\vec{\alpha }}^{\vec{p}} (D^{m} )=\Bigg\{f\in H(D^{m} ):$
\[\left(\int _{D}...\left(\int _{D}\left|f(\omega _{1} ,...,\omega _{m} )\right|^{p_{1} }   \delta ^{\alpha _{1} } (\omega _{1} )d\upsilon (\omega _{1} )\right)^{\frac{p_{2} }{p_{1} } } ...\right)^{\frac{1}{p_{m} } } <+\infty \Bigg\},\] 
for $p_{i+1}\ge p_i,$ $ p_i\ge 1,$ $\alpha>-1,$ $i=1,\dots ,m,$ the following result is valid.

In the following theorem we show that Traces of $A_{\vec{\alpha }}^{\vec{p}} (D^{m} )$ analytic spaces are equal to $A_{\gamma }^{p_{m} } (D)$ spaces.

\begin{theorem} {\rm (see \cite{ShK2015I}).} Let $p_{i+1} \ge p_{i} ,$ $p_{j} \ge 1,j$ $=1,...,m.$ 

Let $\gamma =\alpha _{m} +\sum \limits_{j=1}^{m-1}(n+1+\alpha _{j} )\frac{p_{m} }{p_{j} } , $ $\alpha _{j} >-1,$ $j=1,...,m.$ 

Then $Trace(A_{\vec{\alpha }}^{\vec{p}} (D^{m} ))=A_{\gamma }^{p_{m} } (D),$ $1\le p<+\infty ,$ $j=1,...,m.$
\end{theorem}

\begin{remark}  
{\rm For $D = B$ (unit ball) the theorem 2.5 is known (see \cite{ShM2009A}).}
\end{remark} 

In \cite{ShM2015J} obtained also new sharp estimates for traces in Bergman type spaces of analytic spaces in strongly pseudoconvex domains with smooth boundary.

Let 

$$A_{\vec{\beta }}^{p} (D^{m} )=A^{p} (D^{m} ,\vec{\beta })=\{f\in H(D^{m} ):$$
$$\int _{D}...\int _{D}\left|f(\omega _{1} ,...,\omega _{m} )\right|^{p} \prod_{j=1}^m\delta ^{\beta _{j} } (\omega _{j} )dm(\omega _{1} )...dm(\omega _{m} )<\infty \},$$
where $0<p<\infty ,\beta _{j} >-1,j=\overline{1,m}.$

In the following theorem we provide a sharp result on traces of $A^{p} (D^{m} ,\vec{\beta })$ analytic spaces.

\begin{theorem} {\rm (see \cite{ShM2015J}).} Let $0<p<\infty ,\beta _{j} >-1,j=\overline{1,m}.$ 

Then $Trace\, A^{p} (D^{m} ,\vec{\beta }) = A^{p} (D,\sum \limits _{j=1}^{m}\beta _{j} +(m-1)(n+1)). $
\end{theorem}

We define also the Besov type analytic function spaces on $D^{m} \subset {\mathbb C}^{m} $ as follows 
\[A_{\delta ,\vec{k}}^{p} (D^{m} )=\{ f\in H(D^{m} ):(\sum _{\left|\alpha _{1} \right|\le k_{1} }... \sum _{\left|\alpha _{m} \right|\le k_{m} }\int _{D}...\int _{D}\left|D^{\alpha _{1} ...\alpha _{m} } f\right|   ^{p} \times \] 
\[\times (-\rho _{1} )^{\delta -1} ...(-\rho _{m} )^{\delta -1} d\upsilon (\zeta _{1} ,...,\zeta _{m} ))^{\frac{1}{p} } <+\infty \} ,\] 
$\rho _{j} (z)=\rho (z_{j} ),$ $j=1,...,m,$ $0<p<+\infty ,$ $\delta >0.$

Traces of $A_{t+1,\vec{k}}^{p} (D^{m} )$ classes for $p\le 1$ can be described in a following manner.

\begin{theorem} {\rm (see \cite{ShK2015I}).} $A_{mt+m(n+1)-(n+1)-\sum\limits _{j=1}^{m}k_{j}  }^{p} (D)\subset Trace(A_{t+1,\vec{k}}^{p} (D^{m} )),$ for $p\le 1,$ $t>t_{0} ,$ where $t_{0}$ is large enough.
\end{theorem}

Trace operator in a theorem below is related with certain new interesting sharp embedding theorem in the context of bounded strongly pseudoconvex domains for Bergman type spaces of analytic functions which may have independent interest.

\begin{theorem} {\rm (see \cite{ShK2015I}).} Let $0<p<+\infty,$ $s_j>-1,$ $j=1,\dots ,m.$ Let $\mu $ be a positive Borel measure on $D.$ Then the following are equivalent:
\begin{enumerate}
\item {$\mu (\Delta _{k} )\le C\left|\Delta _{k} \right|^{m+\frac{1}{n+1} \sum \limits_{j=1}^{m}s_{j}  } ,k\ge 1;$}

\item { $Trace(A_{\vec{s}}^{p} )\hookrightarrow L^{p} (D,d\mu );$}

\item { $\int _{D}\prod\limits _{j=1}^{m}\left|f_{j} (z)\right|^{p} d\mu (z)\le c\prod \limits_{j=1}^{m}\left\| f\right\| _{A_{s_{j} }^{p} }^{p} , $ $f_{j} \in A_{s_{j} }^{p} ,$ $1\le j\le m.  $}
\end{enumerate}
\end{theorem}

Next, we define Herz-type spaces in products of pseudoconvex domains $D$ . Note here dealing with weights and measures it is clear from the context we consider products of weights and products of measures, and integration on product domains, though it is not indicated below. Let
\[H_{\alpha}^{p, q}\left(D^{m}\right)=\left\{f \in H\left(D^{m}\right): \sum_{k=1}^{\infty}\left(\int_{B_{0}\left(z_{k}, r\right)}|f(\omega)|^{p}(\delta(\omega))^{\alpha} d V(\omega)\right)^{\frac{q}{p}}<\infty\right\}, \]
where $\delta(z)=\operatorname{dist}(z, \partial D)$ and where $0<p, q<\infty$, $\alpha>-1$ with usual modification when $\alpha$ is vector. Note for $p=q$ we have classical Bergman space and these are Banach spaces when $\min (q, p)>1$ and complete metric spaces for other values of parameters and here $\left\{z^{k}\right\} \subset D$ are certain fixed $r$-lattices in bounded pseudoconvex domains with smooth boundary $D.$
Note obviously these analytic spaces in higher dimension depend on $\left\{z_{n}\right\}$ sequences, but we omit this in names of spaces. We denote here by $d V(z)$ the normalized Lebeques measure on products of bounded pseudoconvex domain with smooth boundary.

In the following theorem we show that Traces of $H_{\vec{\beta}}^{p, q}(D^{m})$ are equal to $A_{\alpha_{1}}^{q}(D)$ spaces.

\begin{theorem} {\rm (see \cite{ShZ2015P}).} Let $F \in H_{\vec{\beta}}^{p, q}\left(D^{m}\right), \quad 0<p$, $q<\infty, n>n_{0}, q \leq p$, and let $\beta_{0}, n_{0}$ be large enough. Then $F(z, \ldots, z) \in A_{a_{1}}^{q}(D)$ that is
\[\int_{D}|F(z, \ldots, z)|^{q} \delta(z)^{\alpha_{1}} d V(z)<+\infty,\]
where $\alpha_{1}=\sum_{j=1}^{m}\left(\beta_{j}+n+1\right) \frac{p}{q}-(n+1).$

Let $f \in A_{\alpha_{1}}^{q}(D)$. Then we can find a function $F$, so that $F \in H_{\vec{\beta}}^{p, q}\left(D^{m}\right)$ and $F(z, \ldots, z)=f(z), z \in D,$ so
\[Trace (H_{\vec{\beta}}^{p, q})(D^{m})=(A_{\alpha_{1}}^{q})(D).\]
\end{theorem}

\section{\large Traces in various spaces of analytic functions in tubular domains}

In this section we consider some interesting problems on traces in various spaces of analytic functions in tubular domains.

Tubular domains are very general unbounded Siegel domains and they were studied in many recent papers of B. Sehba and his various coauthors in recent two decades (see \cite{Sehb2008}, \cite{Sehb2009}). We list these papers in references.
The definition of Traces of analytic function spaces in tube is similar to those given in definition \ref{def1}. In proofs of all our trace theorems below nice properties of $r$-lattices in tube were used.
These important $r$-lattices were invented by A. Bonami and her coauthors (see \cite{Bek2001}). We list many sharp trace theorems below for various Bergman and Herz type function spaces of analytic function. All definitions of these analytic function spaces are given below.

We add the main notations from \cite{Sehb2008}, \cite{Bek2001}, \cite{Bek1990}. 

Let $T_{\Omega} = V + i\Omega$ be the tube domain over an irreducible symmetric cone $\Omega$ in the complexification $V^{\mathbb C}$ of an $n$-dimensional Euclidean space $\tilde V.$ 

We denote the rank of the cone $\Omega$ by $r$ and by $\Delta$ the determinant function on $\tilde V$ . Letting $\tilde V = \mathbb R^n,$ we have as an example of a symmetric cone on $\mathbb {R}^n$ the Lorentz cone $\Lambda_n$ defined for $n \ge 3$ by $\Lambda_n = \{y \in \mathbb {R}^n: y_1^2- \dots -y^2_n > 0, y_1 > 0\}.$ 

Also, we denote $m$-products of tubes by $T^m_{\Omega} = T_{\Omega}\times \dots \times T_{\Omega},$ $\in \mathbb{N}.$ The space of all analytic function
on this new product domain which are analytic by each variable separately will be denoted by $H(T^m(\Omega)).$  

In \cite{ShP2013J} for the first time in literature we consider known problem related with Trace estimates in spaces of analytic functions in unbounded domains in $\mathbb C^n$, namely in tube domains over symmetric cones. In \cite{ShP2015U} authors are interested on properties of certain analytic subspaces of $H\left(T_{\Omega}^{m}\right)$ and obtained new sharp estimates for traces in Bergman type spaces of analytic functions in typical Siegel domains.

By $m$ we denote below a natural number bigger than 1. For $\tau \in \mathbb{R}_{+}$ and the associated determinant function $\Delta(x)$ we set

\[A_{\tau}^{\infty}\left(T_{\Omega}\right)=\left\{F \in \mathcal{H}\left(T_{\Omega}\right):\|F\|_{A_{\tau}^{\infty}}=\sup _{x+i y \in T_{\Omega}}|F(x+i y)| \Delta^{\tau}(y)<\infty\right\}\]
(see \cite{Bek2001} and references there).

For $1 \leq p, q<+\infty$ and $\nu \in \mathbb{R}$, and $\nu>\frac{n}{r}-1$ we denote by $A_{\nu}^{p, q}\left(T_{\Omega}\right)$ the mixed-norm weighted Bergman space consisting of analytic functions $f$ in $T_{\Omega}$ that
\[\|F\|_{A_{\nu}^{p, q}}=\left(\int_{\Omega}\left(\int_{\widetilde{V}}|F(x+i y)|^{p} d x\right)^{q / p} \Delta^{\nu}(y) \frac{d y}{\Delta(y)^{n / r}}\right)^{1 / q}<+\infty.\]

It is known the $A_{\nu}^{p, q}\left(T_{\Omega}\right)$ space is nontrivial if and only if $\nu>\frac{n}{r}-1.$ When $p=q$ we write $A_{\nu}^{p,q}\left(T_{\Omega}\right)=A_{\nu}^{p}\left(T_{\Omega}\right).$

To define related two Bergman-type spaces $A_{\nu}^{p}\left(T_{\Omega}^{m}\right)$ and $A_{\tau}^{\infty}\left(T_{\Omega}^{m}\right)(\nu$ and $\tau$ can be also vectors) in m-products of tube domains $T_{\Omega}^{m}$ we follow standard procedure. For example, we set for all $z_{j}=x_{j}+y_{j}, $ $\tau_{j} \in R,$ $j=1, \ldots, m,$ $ F(z)=F\left(z_{1}, \ldots, z_{m}\right),$ $|F(x+i y)|=\left|F\left(x_{1}+i y_{1}, \ldots, x_{m}+i y_{m}\right)\right|,$ $\tau=\left(\tau_{1}, \ldots, \tau_{m}\right),$
\[A_{\tau}^{\infty}\left(T_{\Omega}^{m}\right)= \left\{F \in H\left(T_{\Omega}^{m}\right):\|F\|_{A_{\tau}^{\infty}}=\sup _{x+i y \in T_{\Omega}^{m}}|F(x+i y)| \Delta^{\tau}(y)<\infty\right\}, \]
where  $\Delta^{\tau}(y)$ is a product of m-onedimensional $\Delta^{\tau_{j}}\left(y_{j}\right)$ functions, $j=1, \ldots, m$. 

Similarly the Bergman space $A_{\tau}^{p}$ can be defined on products of tubes for all $\tau=$ $\left(\tau_{1}, \ldots, \tau_{m}\right), \tau_{j}>\frac{n}{r}-1, j=1, \ldots, m$. It can be shown that all spaces are Banach spaces. Replacing above simply A by L we will get as usual the corresponding larger space of all measurable functions in products of tubes over symmetric cone with the same quazinorm.

As previously in case of analytic functions in unit disk, polydisk, unit ball and upperhalfspace $\mathbb C_+$ and as in proofs of trace theorems in case of spaces of harmonic functions in Euclidean space the role of the Bergman representation formula is crucial. 

In the following theorem we completely describe traces of Bergman spaces.

\begin{theorem} {\rm (see \cite{ShP2013J}).} Let $f \in A_{\nu}^{p}\left(T_{\Omega}^{m}\right), 1 \leqslant p<\infty, \nu \in R^{n}, \nu_{j}>\nu_{0}$ for fixed $\nu_{0}=\nu_{0}(p, n, r, m)$, for all $j=1, \ldots, m$. Then $f(z, \ldots, z) \in A_{s}^{p}$, where $s=\sum_{j=1}^{m}\left(\nu_{j}-\frac{n}{r}\right)+2 \frac{n}{r}(m-1)$ with related estimates for norms. And for all $\frac{n}{r} \leqslant p_{1}$, where $p_{1}$ is a conjugate of $p$ the reverse is also true. For each $g$ function $g \in A_{s}^{p}\left(T_{\Omega}\right)$ there is an $F$ function, $F(z, \ldots, z)=g(z), F \in A_{\nu}^{p}\left(T_{\Omega}^{m}\right)$.
Let in addition
\[\left(T_{\beta} f\right)\left(z_{1}, \ldots, z_{m}\right)=C_{\beta} \int_{T_{\Omega}} f(w) \prod_{j=1}^{m} \Delta^{-t}\left(\left(z_{j}-\bar{w}\right) / i\right) d V_{\beta}(w),\]
$m t=\beta+\frac{n}{r}, z_{j} \in T_{\Omega}, j=1, \ldots, m$.

Then the following asertions hold for all $\beta$ so that $\beta>\beta_{0}$ for some fixed large enough positive number $\beta_{0}$. The $T_{\beta}$ Bergman-type integral operator (expanded Bergman projection) maps $A_{s}^{p}\left(T_{\Omega}\right)$ to $A_{\nu}^{p}\left(T_{\Omega}^{m}\right), \nu=\left(\nu_{1}, \ldots, \nu_{m}\right), \nu_{j}>\nu_{0}, j=1, \ldots, m$.
\end{theorem}

A complete analogue of this theorem is true also for $p=\infty$ case.

In the following theorem we completely describe traces of Bloch
spaces.

\begin{theorem} {\rm (see \cite{ShP2013J}).} Let $f \in A_{\nu}^{\infty}\left(T_{\Omega}^{m}\right), \nu \in R^{n}, \nu_{j}>\frac{n}{r}-1$, for all $j=1, \ldots, m$. Then $f(z, \ldots, z) \in A_{s}^{\infty}$, where $s=\sum_{j=1}^{m} \nu_{j}$. And the reverse is also true - for each $g$ function $g \in A_{s}^{\infty}\left(T_{\Omega}\right)$ there is an $F$ function, $F(z, \ldots, z)=g(z), F \in A_{\nu}^{\infty}\left(T_{\Omega}^{m}\right)$. 

Let in addition
\[\left(T_{\beta} f\right)\left(z_{1}, \ldots, z_{m}\right)=C_{\beta} \int_{T_{\Omega}} f(w) \prod_{j=1}^{m} \Delta^{-t}\left(\left(z_{j}-\bar{w}\right) / i\right) d V_{\beta}(w),\]
$m t=\beta+\frac{n}{r}, z_{j} \in T_{\Omega}, j=1, \ldots, m$. Then the following assertions hold for all $\beta, \beta>\beta_{0}$ for some fixed large enough positive number $\beta_{0}$.

The $T_{\beta}$ Bergman-type integral operator (expanded Bergman projection) maps $A_{s}^{\infty}\left(T_{\Omega}\right)$ to $A_{\nu}^{\infty}\left(T_{\Omega}^{m}\right), \nu=\left(\nu_{1}, \ldots, \nu_{m}\right), \nu_{j}>\frac{n}{r}-1, j=1, \ldots, m, s=\sum_{j=1}^{m} \nu_{j}$.
\end{theorem}

We denote by $A_{\vec{\nu}}^{\vec{p}}\left(T_{\Omega}^m\right)$, $1 \leq p_{j}<+\infty$ and $\nu_{j} \in \mathbb{R},$ $\nu_{j}>\frac{n}{r}-1, j=1, \ldots, m,$ the mixed-norm weighted Bergman space consisting of analytic functions $f$ in $T_{\Omega}$ that
\[
A_{\vec{\nu}}^{\vec{p}}(T_{\Omega}^{m})=\Big\{f \in H(T_{\Omega}^{m}):\Big(\int_{T_{\Omega}} \ldots\Big(\int_{T_{\Omega}}\left|f\left(z_{1}, \ldots, z_{m}\right)\right|^{p_{1}}\left[\Delta\left(y_{1}\right)\right]^{\nu_{1}-n / r} d x_{1} d y_{1}\Big)^{\frac{p_{2}}{p_{1}}} \ldots\]
\[\ldots\left[\Delta\left(y_{m}\right)\right]^{\nu_{m}-n / r} d x_{m} d y_{m}\Big)^{\frac{1}{p_{m}}}<+\infty\Big\}.
\]

This is a Banach space. Replacing above simply A by L we will get as usual the corresponding larger space of all measurable functions in tube over symmetric cone with the same norm.

We define a new Banach space (analytic Herz-type spaces) in tubular domains over symmetric cones as follows. 
Let first 
$$
L_{\tau}^{p}\left(T_{\Omega}^{m}\right)=\left\{f\text{-loc. integrable in } T_{\Omega}:\|f\|_{L_{\tau}^{p}}=\int_{T_{\Omega}}|f(z)|^{p} \prod_{j=1}^{m} \Delta^{\tau_{j}-\frac{n}{r}}\left(y_{j}\right) d x_{j} d y_{j}<\infty,\right\},$$
$\displaystyle \tau_{j}>\frac{n}{r}-1, j=1, \ldots, m$.

To define the next space of functions we remind the reader that the family of Bergman balls $B_{\delta}(z)$ forms an $r$-lattice in tubular domain $T_{\Omega},$ we denote by $B_{\delta}^{m}(z)$ standart $m$-cartesian product of such Bergman balls in $\mathbb C^{m}, $ $z=\left(z_{1}, \ldots, z_{m}\right)$.

By $S_{\nu, \tau}^{p}$ we denote all $f$ functions analytic in $T_{\Omega}^{m}$ so that

$$
\int_{B_{\delta}^{m}(z)}\left|f\left(z_{1}, \ldots, z_{m}\right)\right| \prod_{j=1}^{m} \Delta^{s_{j}}\left(y_{j}\right) d x_{j} d y_{j}
$$
belongs to $L_{\tau_{1}, \ldots, \tau_{m}}^{p}\left(T_{\Omega}^{m}\right),$ where $\displaystyle s_{j}=\nu- \frac{n}{r},$ $\displaystyle \nu_{j}>\frac{n}{r}-1, $ $\displaystyle \tau_{j}>\frac{n}{r}-1, j=1, \ldots, m,$ and $1 \leq p<\infty$ with the same norm (see \cite{Sehb2009}, \cite{Far1994}).

In the following theorem we completely describe traces of general mixed norm spaces. 

\begin{theorem} {\rm (see \cite{ShP2015U}).} Let $f \in A _{\vec{\nu}}^{\vec{p}}\left(T_{\Omega}^{m}\right), $ $ 1 \leq p_{j}<\infty, $ $\nu \in R^{m},$ $ \nu_{j}>\nu_{0}$ for fixed $\nu_{0}=\nu_{0}(n, p, r, m),$ $j=1, \ldots, m$. 

Then $f(z, \ldots, z) \in A_{s}^{p_{m}}, $ $\displaystyle s=\nu_{m}-\frac{n}{r}+\sum_{j=1}^{m-1}\left(\nu_{j}+\frac{n}{r}\right) \frac{p_{m}}{p_{j}}$ with related estimates for norms. And for all $\tilde{p}_{1},$ $\displaystyle \frac{1}{\tilde{p}_{1}}+\frac{1}{p_{m}}=1,$ $ n / r \leq \tilde{p}_{1},$ $j=1, \ldots, m,$ for some fixed $k$ large enough, $\nu_{j}>k$ the reverse is also true. For each function $g \in A_{s}^{p_{m}}\left(T_{\Omega}\right)$ there is an $F$ function $F(z, \ldots, z)=g(z), F \in A_{\vec{\nu}}^{\vec{p}}\left(T_{\Omega}^{m}\right)$.

Let in addition
$$T_{\beta}(f)\left(z_{1}, \ldots, z_{m}\right)=C_{\beta} \int_{T_{\Omega}} f(w)\left(\prod_{j=1}^{m} \Delta^{-t}\left(\frac{z_{j}-\bar{w}}{i}\right)\right) d V_{\beta}(w),$$ 
$\displaystyle mt=\beta+\frac{n}{r},$ $ z_{j} \in T_{\Omega},$ $ j=1, \ldots, m$.

Then following assertion hold for all $\beta$, so that $\beta>\beta_{0}$ for some fixed large enough positive number $\beta_{0}$. The $T_{\beta}$ Bergman-type integral operator (expanded Bergman projection) maps $A_{s}^{p_{m}}\left(T_{\Omega}\right)$ to $A_{\vec{\nu}}^{\vec{p}}\left(T_{\Omega}^{m}\right),$ $ \nu=\left(\nu_{1}, \ldots, \nu_{m}\right),$ $ \nu_{j}>\nu_{0},$ $ j=1, \ldots, m$.
\end{theorem}

For $p_{j}=p, j=1, \ldots, m$, case where $1 \leq p<\infty$ this theorem was proved in \cite{ShP2013J}. For $m=1$ this theorem is obvious. For $n=1$ we have upperhalfspace $\mathbb C_{+}$, this theorem was proved in \cite{ShM2010S}. 

In the theorems \ref{th3.4}, \ref{th3.5} we completely describe traces of new Herz  spaces. 

\begin{theorem} \label{th3.4} {\rm (see \cite{ShP2015U}).} Let $\nu_{j}>\nu_{0}$ and $\tau_{j}>\tau_{0}$ for some fixed positive numbers $\nu_{0}=\nu_{0}(p, n, r, m)$ and $\tau_{0}=\tau_{0}(p, n, r, m),$ $ 1 \leq p<\infty$. Let $f \in S_{\nu, \tau}^{p}\left(T_{\Omega}^{m}\right).$ 

Then $f(z, \ldots, z)$ belongs to $A_{s}^{p}\left(T_{\Omega}\right)$ where $\displaystyle s=\sum_{j=1}^{m}\left(\nu_{j}+\frac{n}{r}\right) p+\sum_{j=1}^{m}\left(\tau_{j}+\frac{n}{r}\right)-\frac{n}{r}$ and for every $f$ function $f \in A_{s}^{p}\left(T_{\Omega}\right)$ there is an $F$ function $F \in S_{\nu, \tau}^{p}$ so that $F(z, \ldots, z)=f(z)$ for all $\displaystyle \frac{n}{r} \leq p^{\prime},$ $\displaystyle \frac{1}{p}+\frac{1}{p^{\prime}}=1$.

Let in addition
$$
\left(T_{\beta} f\right)\left(z_{1}, \ldots, z_{m}\right)=C_{\beta} \int_{T_{\Omega}} f(w) \prod_{j=1}^{m} \Delta^{-t}\left(\frac{z_{j}-\bar{w}}{i}\right) d V_{\beta}(w),
$$
$\displaystyle mt=\beta+\frac{n}{r},$ $ z_{j} \in T_{\Omega}, $ $j=1, \ldots, m$. 

Then the following assertions holds for all $\beta$ so that $\beta>\beta_{0}$ for some fixed large enough positive number $\beta_{0}$. The $T_{\beta}$ Bergman projection-type integral operator maps $A_{s}^{p}\left(T_{\Omega}\right)$ to $S_{\nu, \tau}^{p}\left(T_{\Omega}^{m}\right),$ $ \nu=\left(\nu_{1}, \ldots, \nu_{m}\right),$ $ \tau=$ $\left(\tau_{1}, \ldots, \tau_{m}\right),$ $ \nu_{j}>\nu_{0},$ $ \tau_{j}>\tau_{0},$ $ j=1, \ldots, m$.
\end{theorem}

For $p=1$ this theorem was proved in \cite{ShP2013J} in tubular domains.

Further, we define Herz-type spaces in products of tubular domains $T_{\Omega}$. In \cite{ShZ2015P} the authors provided new sharp results on traces in some new analytic Herz-type spaces in some general domains (unbounded typical Siegel domains of the second type) in $\mathbb {C}^n$ (tubular domains). These spaces serve as very natural generalizations of classical Bergman spaces. Note here dealing with weights and measures it is clear from the context we consider products of weights and products of measures. Let

\[H_{\alpha}^{p, q}\left(T_{\Omega}^{m}\right)=\left\{f \in H\left(T_{\Omega}^{m}\right): \sum_{k=1}^{\infty}\left(\int_{B_{r}\left(z_{k}, r\right)}|f(\omega)|^{p}(\operatorname{Im}^{\alpha} \omega) d V(\omega)\right)^{\frac{q}{p}}<\infty\right\},\]
where $0<p, q<\infty$, $\alpha>-1$ with usual modification when $\alpha$ is vector. Note for $p=q$ we have classical Bergman space and these are Banach spaces when $\min (q, p)>1$ and complete metric spaces for other values of parameters and here $\left\{z^{k}\right\} \subset T_{\Omega}$ are certain fixed $r$-lattices in tubular domains over symmetric cones $T_{\Omega}$ which were discussed above. Note obviously these analytic spaces in higher dimension depend on $\left\{z_{n}\right\}$ sequences, but we omit this in names of spaces. We denote here by $d V(z)$ the normalized Lebeques measure on products of tubular domains.

In the following theorem the authors provided a sharp result on traces of $H_{\vec{\beta}}^{p, q}\left(T_{\Omega}^{m}\right)$ analytic spaces.

\begin{theorem} \label{th3.5} {\rm (see \cite{ShZ2015P}).} Let $F \in H_{\vec{\beta}}^{p, q}\left(T_{\Omega}^{m}\right), 1<p,q<\infty$, $p \geq q, \beta_{0}, \alpha_{1}$ be large enough. Then $F(z, . ., z) \in A_{\alpha_{1}}^{q}\left(T_{\Omega}\right)$ that is
\[\int_{T_{\Omega}}|F(z, \ldots, z)|^{q}(\operatorname{Imz})^{\alpha_{1}} d V(z)<\infty,\]
where $
\alpha_{1}=\sum_{j=1}^{m}\left(\beta_{j}+\frac{2 n}{r}\right) \frac{p}{q}-2 \frac{n}{r}.$

Let $f \in A_{\alpha_{1}}^{q}\left(T_{\Omega}\right)$. Then we can find a $F$ function, so that $F \in H_{\vec{\beta}}^{p, q}\left(T_{\Omega}^{m}\right)$ and $F(z, \ldots, z)=f(z)$, $z \in T_{\Omega}$, so
$$
Trace (H_{\vec{\beta}}^{p, q}(T_{\Omega}^{m}))=A_{\alpha_{1}}^{q}(T_{\Omega}).
$$
\end{theorem}

The \cite{ShM2019R} presents some results for unbounded tubular domains over symmetric cones. Let
\[A_{\vec{s}}^{p}=\left\{f \in H\left(T_{\Omega}^{m}\right): \int_{T_{\Omega}} \cdots \int_{T_{\Omega}}\left|f\left(z_{1}, \ldots, z_{m}\right)\right|^{p} \prod_{j=1}^{m} \Delta^{s_{j}-\frac{n}{r}}(\operatorname{Im} z) d V(z)<+\infty\right\},\]
where $s=\left(s_{1}, \ldots, s_{m}\right), 0<p<\infty, s_{j}>\frac{n}{r}-1, j=1, \ldots, m$. Let also further $L^{p}\left(T_{\Omega}, d \mu\right)$ be a space of all measurable functions on $T_{\Omega}$ so that $\int_{T_{\Omega}}|f(z)|^{p} d \mu(z)<\infty,$ for any positive Borel measure $\mu$ on $T_{\Omega}$.

In the following theorem certain new interesting sharp embedding theorems related with the trace operator are formulated.
Similar type sharp embedding theorems related with the trace operator are valid also in other domains. They were proved by first author in recent years.

\begin{theorem} {\rm (see \cite{ShM2019R}).} Let $0<p<+\infty,$ $s_j>-1,$ $j=1,\dots ,m.$ Let $\mu $ be a positive Borel measure on $T_{\Omega } $. Then the following are equivalent:
\begin{enumerate}
\item { $\mu (\Delta _{k} )\le C\left|\Delta _{k} \right|^{m+\left(\frac{2n}{r} \right)^{-1} \sum \limits_{j=1}^{m}s_{j}  } ,k\ge 1;$}

\item { $Trace(A_{\vec{s}}^{p} )\hookrightarrow L^{p} (T_{\Omega } ,d\mu );$}

\item {$\int _{T_{\Omega } }\prod \limits_{j=1}^{m}\left|f_{j} (z)\right|^{p} d\mu (z)\le c\prod\limits _{j=1}^{m}\left\| f\right\| _{A_{s_{j} }^{p} }^{p} ,$ $f_{j} \in A_{s_{j} }^{p} ,$ $1\le j\le m. $} 
\end{enumerate}
\end{theorem}

In \cite{ShM2022R} gives new exact theorems on traces of analytic classes of BMOA type in tubular regions over symmetric cones under one additional condition on the so-called Bergman kernel in these domains. An essential role is played by the so-called $r$ -lattices for tubular regions obtained earlier in the works of foreign authors. The theorems of this paper are complete analogs of recent exact results on the traces of classes of BMOA, obtained earlier by the first author in polyballs and in bounded strongly pseudoconvex domains.

Let $T_{\Omega}$ be a tubular domain over symmetric cones with smooth boundary in $\mathbb C^n.$
Let 
\[M^p_{r_1 ,...,r_m,\tau,s_1 ,...,s_m } (T_{\Omega}^m )=\{ f \in H(T_{\Omega}^m):\mathop{\sup }\limits_{\omega \in T_{\Omega} } \delta^{m\tau } (\omega )\int _{T_{\Omega} }...\int _{T_{\Omega} }\left|f(z_1 ,...,z_m)\right|^p  \times  \]
\[\times \prod \limits_{j=1}^{m}\left|B_{s_j } (z_j ,\omega )\right| \prod \limits_{j=1}^{m}\delta ^{s_j } (z_j )dV(z_j)<+\infty \} , \] 
$s_{j} =\tau +r_{j} ,$ $j=1,...,m,$ where $\tau >0,$ $s_{j} >-1,$ $r_{j} \ge 0,$ $j=1,...,m,$ $0<p<+\infty .$

These are Banach spaces for all $p \geq 1$ and complete metric spaces for other values of $p$. The following theorems for unit ball case can be seen in \cite{ShM2010S}. 

In the theorems \ref{th3.7}, \ref{th3.8} we provide a sharp result on traces of $M^p_{ r_1,...,r_m,\tau ,s1,...,sm} (T_{\Omega}^m)$ spaces.

These two following sharp trace results are valid under certain composition formula in tubular domain over simmetric cones which can be seen in mentioned papers below and which is valid for unit polydisk and unit ball and bounded strongly pseudoconvex domains.

\begin{theorem} \label{th3.7} \rm (see \cite{ShM2022R}). Let $p > 1,$ $ \tau \in (0,+\infty ),$ $ r_j \in \mathbb N ,$ $r =\sum\limits_{j=1}^m {r_j },$ $ s_j > - 1,$ $ j = 1, . . . , m,$ $t = (m - 1)(2n/r)+ \sum\limits_ {j=1}^m {s_j}.$  
Then $Trace(M^p_{ r_1,...,r_m,\tau ,s1,...,sm} (T_{\Omega}^m)) = M^p_{r,t,\tau m}(T_{\Omega})$ for all $n,n/r>n_0,$ where $n_0=n_0(p,\tau,r_1,…,r_m,m).$
\end{theorem}

\begin{theorem} \label{th3.8} {\rm (see \cite{ShM2022R}).} Let $p \leq 1,$ $ \tau \in (0,+\infty ),$ $ r_j \in \mathbb N ,$ $r =\sum\limits_{j=1}^m {r_j },$ $ s_j > - 1,$ $ j = 1, . . . , m,$ $\displaystyle \frac{r_j}{p} \in \mathbb N, j=1,…m, $ $t = (m - 1)(2n/r)+ \sum\limits_ {j=1}^m {s_j}.$  
Then $Trace(M^p_{ r_1,...,r_m,\tau ,s1,...,sm} (T_{\Omega}^m)) = M^p_{r,t,\tau m}(T_{\Omega})$ for all $n,n/r>n_0,$ where $n_0=n_0(p,\tau,r_1,…,r_m,m).$
\end{theorem}

All sharp trace theorems of this forth section are valid with rather similar proofs in analytic spaces of several variables in various domains and they were provided earlier by first author recently in various papers based on nice properties of so called $r$-lattices in these domains in $\mathbb C^n.$

\section{\large Traces in various spaces of harmonic functions in products of upper half spaces}

In this section we consider some interesting problems on traces in various spaces of harmonic functions in products of upper half spaces.

Let $B$ be the unit ball in ${\mathbb R}^{n} $, that is 
\[B=\left\{x=(x_{1} ,...,x_{n} )\in {\mathbb R}^{n} :\left|x\right|=\left(\sum \limits_{i=1}^{n}\left|x_{i} \right|^{2}  \right)^{\frac{1}{2} } \le 1\right\},\] and $S^{n-1} =\partial B=\left\{x=(x_{1} ,...,x_{n} )\in {\mathbb R}^{n} :\left|x\right|=1\right\},$ and $x=rx',r=|x|\in (0,1),x'\in S^{n-1} .$

We consider Banach Bergman spaces in $B$ (see, for example, \cite{DjShF1988}) 
\[A_{\alpha }^{p} (B)=\left\{f\in h(B):\left\| f\right\| _{p,\alpha } :=\left(\int _{0}^{1}\int _{S^{n-1} }\left|f(rx')\right|^{p} (1-r)^{\alpha } r^{n-1} drdx'  \right)^{\frac{1}{p} } <+\infty \right\},\] 
when $1\le p<+\infty $, $0\le \alpha <+\infty $, and $p=+\infty $, $0\le \alpha <+\infty $ we have
\[\left\| f\right\| _{\infty ,\alpha } :=\mathop{\sup }\limits_{x\in B} \left|f(x)\right|\left(1-\left|x\right|\right)^{\alpha } <+\infty .\] 

Let also $B^m=B\times...\times B.$  We denote $A_{\alpha _{1} ,...,\alpha _{m} }^{\infty } (B^m),$  $\alpha _{j} >0,$ $j=1,...,m,$ the space of all functions $f$ harmonic by each variable in $B$ such that 
\[\mathop{\sup }\limits_{x_{i} \in B,i=1,...,m} \left|f(x_{1} ,...,x_{m} )\right|\prod _{i=1}^{m}\left(1-\left|x_{i} \right|^{2} \right)^{\alpha _{i} }  <+\infty .\] 

Let $\alpha _{j} >0,$ $j=1,...,m.$ Some information about $Trace\left(A_{\alpha _{1} ,...,\alpha _{m} }^{\infty } (B^{m} )\right)$ can be seen in \cite{ShM2011F}.

Let 
\[hA_{\vec{\alpha }}^{p}  (B^m)=\{ f\in h(B^{m} ):\int _{B}\int _{B}...  \int _{B}\int _{B}\left|f(x_{1} ,...,x_{m} )\right|^{p} \times  \] 
\[\times\prod _{k=1}^{m}\left(1-\left|x_{k} \right|\right)^{\alpha _{k} }  (\left|x_{1} \right|...\left|x_{m} \right|)^{n-1} d\left|x_{1} \right|...d\left|x_{m} \right|dx'_{1} ...dx'_{m} <+\infty \} .\] 

Let $p\ge 1,$ $m\in {\mathbb N},$ $\alpha _{j} >-1,$ $j=1,...,m.$ Some information about $Trace(hA_{\vec{\alpha }}^{p} (B^{m} ))$ can be seen in \cite{ShM2011F}.

We introduce harmonic Bergman classes on ${\mathbb R}^{n+1} $ and on product spaces as follows. 

We set $\mathbb R^{n+1}_+ = \{(x, t) : x \in \mathbb R^n, t > 0\} \subset \mathbb R^{n+1}.$ We denote the points in $\mathbb R^{n+1}_+$ usually by $z = (x, t)$ or $w = (y, s).$ The Lebegue measure is denoted by $dm(z) = dz = dxdt$ or $dm(w) = dw = dyds.$ We also use standard weighted measures $dm_{\lambda}(z) = t^{\lambda} dxdt,$ $\lambda \in \mathbb R.$ 

Let $0<p<+\infty ,$ $\alpha >-1,$
\[A_{\alpha }^{p} \left({\mathbb R}_{+}^{n+1} \right)=\{ f\in h\left({\mathbb R}_{+}^{n+1} \right):\left\| f\right\| _{p,\alpha }^{p} =\int _{{\mathbb R}^{n} }\int _{0}^{\infty }\left|f(x',x_{n+1} )\right|^{p} x_{n+1}^{\alpha }   dx'dx_{n+1} <+\infty \} ,\] 
and 
\[A_{\alpha_1,...\alpha_m}^{p} \left({\mathbb R}_{+}^{n+1} \times ...\times {\mathbb R}_{+}^{n+1} \right)= A_{\vec{\alpha} }^{p} \left(({\mathbb R}_{+}^{n+1})^m\right)=\{ f\in h\left( ({\mathbb R}_{+}^{n+1})^m \right):\left\| f\right\| _{p,\alpha }^{p} =\]
\[=\int _{{\mathbb R}^{n} }...\int _{{\mathbb R}^{n} }\int _{0}^{\infty }...   \int _{0}^{\infty }\left|f(x'_{1} ,...,x'_{m} ,x_{n+1}^{1} ,...,x_{n+1}^{m} )\right|^{p} \times  \] 
\[\times \prod _{k=1}^{m}\left(x_{n+1}^{k} \right)^{\alpha _{k} } dx'_{1} ...dx'_{m}  dx_{n+1}^{1} ...dx_{n+1}^{m} <+\infty \}.\]

Let $0<p<+\infty ,$ $m\in {\mathbb N},$ $\alpha _{j} >-1,$ $j=1,...,m.$ Some information about $Trace(A_{\alpha _{1} ,...,\alpha _{m} }^{p} ({\mathbb R}_{+}^{n+1} \times ...\times {\mathbb R}_{+}^{n+1} ))$ can be seen in \cite{ShM2011F}.

Trace theorems of various type and various estimates of traces of various harmonic function spaces can be seen in \cite{ArSh2011K}, \cite{ArSh2013K}. We formulate some results from these papers.

We set $\mathbb H=\{ (x,t):x\in {\mathbb R}^{n} ,t>0\} \subseteq {\mathbb R}^{n+1} .$ 
Weighted harmonic Bergman spaces on $\mathbb H$ are defined, for $0<p<+\infty$ and $\lambda>-1,$ as usual
\[A_{\lambda }^{p} =A_{\lambda }^{p} (\mathbb H)=\left\{ f\in h(\mathbb H):\left\| f\right\| _{A_{\lambda }^{p} } =\left(\int _{\mathbb H}\left|f(z)\right|^{p}  dm_{\lambda } (z)\right)^{\frac{1}{p}} <+\infty \right\} .\] 

For $\vec{\alpha }=(\alpha _{1} ,...,\alpha _{m} )\in {\mathbb R}^{m} $ we have a product measure $dm_{\vec{\alpha }} $ on $\mathbb H^m $ defined by $dm_{\vec{\alpha }} (z_1 ,...,z_m )=dm_{\alpha _1 } (z_1 )...dm_{\alpha _m } (z_m )$ and we set $L_{\vec{\alpha }}^{p} =L^{p} (\mathbb H^{m} ,dm_{\vec{\alpha }} ),$  $A_{\vec{\alpha }}^{p} =L_{\vec{\alpha }}^{p} \cap h(\mathbb H^{m} ),$ $0<p<+\infty .$ We denote by $\tilde A_{\vec{\alpha }}^{p} $ the subspace of $A_{\vec{\alpha }}^{p} $ consisting of functions which are harmonic in each of the variables $z_1,\dots, z_m$ separately. 

For a function $f: \mathbb H^{m} \to {\mathbb C}$ we define $Tr\, f:{\mathbb H} \to {\mathbb C}$ by $Tr\,f(z)=f(z,...,z).$

The following results are true.

\begin{theorem} {\rm (see \cite{ArSh2011K}).} Let $1<p<+\infty,$ $s_1,\dots ,s_m>-1$ and set $\lambda =(m-1)(n+1)+\sum\limits _{j=1}^{m}s_{j}.$ Then $A_{\lambda }^{p} \subset Trace(\tilde{A}_{\vec{s}}^{p} )\subset Trace (A_{\vec{s}}^{p} )\subset L^{p} ({\mathbb R}_{+}^{n+1} ,dm_{\lambda } ).$

In particular, if $f\in A_{\vec{s}}^{p} $ and if $Tr \,f$ is harmonic, then $Tr\,f \in A^p_{\lambda}.$
\end{theorem}

For $p=\infty ,$ $s_{j} >0,$ $1\le j\le m,$ we define $A_{\vec{s}}^{\infty } $ as the space of all $f\in h(({\mathbb R}_{+}^{n+1} )^{m} )$ such that
\[\left\| f\right\| _{A_{\vec{s}}^{\infty } }=\mathop{\sup }\limits_{z_{j} \in {\mathbb R}_{+}^{n+1} } \left|f(z_{1} ,...,z_{m} )\right|t_{1}^{s_{1} } ...t_{m}^{s_{m} } <+\infty ,\] 
for $m=1,$ $s_{1} =\alpha,$ we use simpler notation $A_{\alpha }^{\infty } $.  

The space of all functions $f(z_1, \dots, , z_m)$ on $(R^{
n+1}_+ )^m$ which are harmonic in each of variables
$z_1,\dots, z_m$ is denoted by $\tilde h((R^{n+1}_
+)^m).$ Let also if $X$ is a space of function harmonic on $\left({\mathbb R}_{+}^{n+1} \right)^{m}, $ then we set $\tilde{h}X=X\cap \tilde{h}(({\mathbb R}_{+}^{n+1} )^{m} ).$

Let $X\subset h(({\mathbb R}_{+}^{n+1} )^{m} ).$ The trace of $X$ is $TraceX=\{ Tr\,f:f\in X\} .$

\begin{theorem} {\rm (see \cite{ArSh2013K}).} Let $s_1,\dots, s_m>-1$ and $\lambda =\sum \limits_{j=1}^{m}s_{j}  .$  Then $A_{\lambda }^{\infty } \subset Trace(\tilde{h}A_{\vec{s}}^{\infty } )$.

Conversely, if $f\in A_{\vec{s}}^{\infty } $ and $Tr\,f$ is harmonic, then $Tr\,f\in A^{\infty}_{\lambda}.$
\end{theorem}

We below formulate a complete analogue of theorem above in case of harmonic function spaces in the unit ball of $\mathbb R^n$.
The formulation is very similar. We refer for definition of objects in this theorem to mentioned paper of the author with M. Arsenovic.

\begin{theorem} {\rm (see \cite{ArSh2013K}).} Let $0<p<1,$ $s_1,\dots ,s_m>-1$ and $\lambda =(m-1)(n+1)+\sum\limits _{j=1}^{m}s_{j}  .$  Then $A_{\lambda }^{p} \subset Trace(\tilde{h}A_{\vec{s}}^{p} )\subset Trace(A_{\vec{s}}^{p} )\subset L^{p} ({\mathbb R}_{+}^{n+1} ,dm_{\lambda } ).$ 

In particular, if $f\in A_{\vec{s}}^{p} $ and if $Tr \,f$ is harmonic, then $Tr\,f \in A^p_{\lambda}.$
\end{theorem}

\section{\large Final remarks and some lines of development and traces of some new BMOA type, Bloch type, Herz type spaces and some new related results}

In this section we collect many interesting new problems questions and various comments concerning trace problem in analytic function spaces of several complex variables not only inbounded strongly pseudoconvex and tubular domains but also in polydisk and unit ball and polyball. This section may be also viewed as a continuation  of the first part of this expository paper in some sence. So we refer the interested reader to that expository paper of the authors concerning various standard objects and definitions of those objects related with the polydisk and polyball.

We define in this section not only various analytic function spaces in various domains in $\mathbb{C}^n$ and give some information on traces of these spaces but also consider maps which are close to trace map and provide some interesting information concerning their action in such type function spaces, obtained previously by first author.

In certain analytic function spaces in product domains below we see well known $\mathcal{D}^\alpha$ differentiation operators in quazinorms and give estimates of traces of such spaces also.

We discuss now some new interesting open questions Traces of mixed norm spaces with quasinorms in the unit polydisk 

\[\int _{0}^{1}...\int _{0}^{1}\left(\int _{T^{n} }\left|f\left(r\zeta \right)\right|^{p} d\zeta  \right)  ^{\frac{q}{p}} \left(1-r\right)^{\alpha } dr.\] 

Even with more general $\omega \left(r\right)$ weights than $\left(1-r\right)^{\alpha } $ were completely described in recent papers \cite{JPSh2009P} for all $0<p,q<+\infty ,\alpha >-1.$

The next very natural question is to completely describe traces of similar type function analytic spaces with the finite quazinorms defined on $B\times ...\times B$ polyball (products of balls)
\[\int _{0}^{1}...\int _{0}^{1}\left(\int _{S}... \int _{S}\left|f\left(r_{1} \zeta _{1} ,...,r_{n} \zeta _{n} \right)\right|^{p} d\zeta  \right)  ^{\frac{q}{p} } \prod _{j=1}^{n}\left(1-r\right)^{\alpha _{j} }  d\vec{r},\]
$S\subset \mathbb C^{m} ,\alpha _{j} >-1,$ $j=1,...,m,$ $0<p,q<+\infty .$ 

Such spaces for $n=1$ are well studied. Note that similar type even more general problems can be also posed in complete analogues of these analytic function spaces but in general tubular domains over symmetric cones $T_{\Omega } $, or in bounded strongly pseudoconvex domains with smooth boundary $D$, namely on products of these domains.

Traces of $H^{p} $ Hardy spaces for all values of $p$ were given in \cite{ShM2008R}, \cite{ShK2015I}, \cite{ShK2016}. The natural open question to find description of traces of mixed norm $H^{\vec{p}} $.

Mixed norm Hardy spaces defined as follows with the help of the following quazinorms in the unit polydisk in $\mathbb C^{n}$
\[\mathop{\sup }\limits_{r_{m} <1} ...\mathop{\sup }\limits_{r_{1} <1} \left(\int _{T}...\left(\int _{T}\left|f\left(\vec{r}\vec{\zeta }\right)\right|^{p_{1} } d\zeta _{1}  \right)^{\frac{p_2}{p_1} } ...d\zeta _{m}  \right)^{\frac{1}{p_m}} ,0<p<+\infty ,j=1,...,m,\] 
(and also weigthed important $H_{\omega }^{\vec{p}} $  generalization).

Some sharp estimates for traces of new these $H^{\vec{p}} $ classes can be seen in \cite{ShK2015I}, \cite{ShK2016}, in \cite{ShM2009J}, \cite{ShM2009A} for some special cases of these spaces in the polyball, in both cases estimates are given for m=2 and particular values of $p_{j} $ in recent paper \cite{Chas2003}.

Note that very similar question may be raised in $H^{\vec{p}} $ type spaces of analytic functions in $V^{n} $ polydisk with quazinorms $\Gamma _{t} \left(\zeta \right)$ is a Lusin cone: $\Gamma _{t} \left(\zeta \right)=\left\{z\in V:\left|1-\bar{\zeta }z\right|<t\left(1-\left|z\right|\right)\right\},t>1,$ (or replacing $\Gamma _{t} \left(\zeta _{j} \right)$ by $r_{j} <1$).
\[\left(\int _{T}\mathop{\sup }\limits_{z_{m} \in \Gamma _{t} \left(\zeta \right)} ...\left(\int _{T}\mathop{\sup }\limits_{z_{1} \in \Gamma _{t} \left(\zeta \right)} \left|f\left(z_{1} ,...,z_{m} \right)\right|^{p_{1} } d\zeta _{1}  \right)^{\frac{p_2}{p_1}} ...d\zeta _{m}  \right)^{\frac{1}{p_m} },\] 
$0<p_{j} <+\infty ,j=1,...,m,$
or
\[\left(\int _{T}\left(\int _{T}\mathop{\sup }\limits_{z_{1} \in \Gamma _{t} \left(\zeta \right)} ...\mathop{\sup }\limits_{z_{m} \in \Gamma _{t} \left(\zeta \right)} \left|f\left(z_{1} ,...,z_{m} \right)\right|^{p_{1} } d\zeta _{1}  \right)^{\frac{p_2}{p_1} } ...d\zeta _{m}  \right)^{\frac{1}{p_m} },\]
$0<p_{j} <+\infty ,j=1,...,m.$

Note many results are valid also for spaces of pluriharmonic functions. We now introduce new mixed norm spaces in the polyballs $B\times\dots \times B$ as follows

\[A_{\alpha _{1} ...\alpha _{m} }^{p_{1} ...p_{m} } \left(B_{n}^{m} \right)=\left\{f\in H(B_{n} \times ...\times B_{n} ):\right. \left\| f\right\| _{A_{\vec{\alpha }}^{\vec{p}} } =\]
\[\left(\int _{B_{n} }\left(1-\left|z_{m} \right|\right)^{\alpha _{m} }   ...\left(\int _{B_{n} }\left|f\left(z_{1} ,...,z_{m} \right)\right|^{p_{1} } \left(1-\left|z_{1} \right|\right)^{\alpha _{1} } d\nu \left(z_{1} \right) \right)^{\frac{p_2}{p_1} } ... d\nu \left(z_{m} \right)\right)^{\frac{1}{p_m}} <+\infty \},\] 
$0<p<+\infty ,\alpha _{j} >-1,j=1,...,m.$

These analytic spaces can be introduced also for $p_{j} =+\infty $. For $1\le p_{j} <+\infty $ we have Banach spaces for other values complete metric spaces.

We give complete description of mixed norm Bergman spaces in the polyball and upperhalfplane below. It will be interesting to provide sharp extensions of these results to general  weighted spaces and to tube and pseudoconvex domains.

\begin{theorem} {\rm (see \cite{ShM2008R}, \cite{ShM2009A}, \cite{ShM2010J}, \cite{ShK2015P}).} Let  $\gamma=\alpha _{m}+\sum _{j=1}^{m-1}\left(n+1+\alpha _{j} \right)\left(\frac{p_{m} }{p_{j} } \right),$ $\alpha _{j} >-1,$ $1<p_{j} <+\infty,$ $j=1,...,m.$ 

We have a sharp result $Trace\left(A_{\alpha _{1} ,...,\alpha _{m} }^{p_{1} ,...,p_{m} } \left(B_{n}^{m} \right)\right)=A_{\gamma}^{p_{m} } \left(B_{n} \right).$ 
\end{theorem}

Let further 
\[A_{\alpha }^{p} \left(\mathbb C_{+}^{m} \right)=\{ f\in H\left(\mathbb C_{+}^{m} \right):\]
\[\left\| f\right\| _{A_{\alpha }^{p}} =\]
\[=\left(\int _{0}^{+\infty }\int _{0}^{+\infty }...  \int _{0}^{+\infty }\int _{-\infty }^{+\infty }\left|f\left(x_{1} +iy_{1} ,...,x_{m} +iy_{m} \right)\right|^{p} \prod _{j=1}^{m}y_{j}^{\alpha _{j} } dx_{j} dy_{j}    \right)^{\frac{1}{p} } <+\infty \} ,\] 
$0<p<+\infty,$ $\alpha _{j} >-1,$ $j=1,...,m.$

\begin{theorem} {\rm (see \cite{ShM2008R}, \cite{ShM2009A}, \cite{ShM2010J}, \cite{ShK2015P}).} Let $0<p<+\infty , $ $j=\sum _{j=1}^{m}\alpha _{j} +2m-2,$ $ \alpha _{j} >-1,$ $j=1,...,m.$

We have $Trace\left(A_{\alpha }^{p} \left(\mathbb C_{+}^{m} \right)\right)=A_{j}^{p} \left(\mathbb C_{+} \right).$
\end{theorem}

We define certain new weighted Hardy spaces in the polyball and Bergman-Sobolev type spaces in the polyball.and extend some estimates and sharp results provided earlier in polydisk previously by other authors. Note here that we refer the reader for well-known  fractional $\mathcal {D}^\alpha$ derivative which appear below to related papers which we cited below. To extend these results to very complicated tubular domains is an interesting and probably rather difficult open problem.

Some nice generalization of Diagonal map can be seen in a paper of D. Clark \cite{Cl1988}.

Let 
\[H_{\alpha ,\beta }^{p} =\left\{f\in H(B^{m} ):\mathop{\sup }\limits_{r<1} \tilde{M}_{p} (\mathcal D^{\beta } f,r)\left(1-r\right)^{\alpha } <+\infty ,\alpha \ge 0,0<p<+\infty ,\beta \in {\mathbb R}\right\},\] 
where $\mathcal D^{\beta } $ is a fractional derivative of $f$, $\mathcal D^{\beta } =\mathcal D_{z_{1} }^{\beta } ...\mathcal D_{z_{m} }^{\beta } f,$

\[\tilde{M}_{p}^{p} \left(f,r\right)=\int _{S}...\int _{S}\left|f(r\zeta _{1} ,...,r\zeta _{m} )\right|^{p}   d\sigma (\zeta _{1} )...d\sigma (\zeta _{m} ),r\in \left(0;1\right).\] 

We define $H^{p} (B^{m} )=H_{0,0}^{p} $ for $0<p+\infty .$

\[A_{t,\vec{\alpha }}^{p} =\left\{f\in H(B^{m} ):\int _{B}...\int _{B}\left|\mathcal D^{\alpha _{1} } ...\mathcal D^{\alpha _{m} } f\right|^{p} \left(1-\left|\tilde{z}\right|\right)^{t} d\nu (z)  <+\infty \right\},\] 
$\alpha _{j} \in {\mathbb R},$ $j=1,...,m,$ $t>-1,$ $0<p<+\infty,$ $\displaystyle \left(1-\left|\tilde{z}\right|\right)=\prod\limits_{k=1}^{m}\left(1-\left|z_{k} \right|\right) .$

\begin{theorem} {\rm (see \cite{ShM2008R}, \cite{ShM2009A}, \cite{ShM2010J}, \cite{ShK2015P}).}

\begin{enumerate}
\item  $Trace\, (H^{P} )=A_{n(m-1)-1}^{p},$ $m,n\in {\mathbb N},$ $0<p<+\infty .$

\item  Let $p\le 1,$ $\alpha \ge \beta \ge 0,$ then
\[A_{\tau }^{p} \subset Trace\, \left(H_{\alpha ,\beta }^{p} \right),\tau =nm+\left(\alpha -\beta \right)pm-\left(n+1\right).\] 
\item  Let $p\le 1,$ $\alpha _{j} \ge 0,$ $\displaystyle t\ge \frac{\sum \limits_{j=1}^{m}\alpha _{j}  }{m} $, then
\[A_{tm+m(n+1)-(n+1)-\sum\limits _{j=1}^{m}\alpha _{j}  }^{p} \subset Trace\, (A_{t,\vec{\alpha }}^{p} ).\] 
\end{enumerate}
\end{theorem}

It will be interesting to extend this result to tube and pseudoconvex domains.

Traces of new analytic function spaces in the polydisk with the following mixed quasinorm
\[\left(\int _{U}...\left(\int _{U}\left|f\left(z_{1} ,...,z_{m} \right)\right|^{p_{1} } \left(1-r_{1} \right)^{\alpha _{1} } dr_{1} d\zeta _{1}  \right)^{\frac{p_2}{p_1} } ...\left(1-r_{m} \right)^{\alpha _{m} } dr_{m} d\zeta _{m}  \right)^{\frac{1}{p_m} },\]
$0<p_{j} <+\infty ,\alpha _{j} >-1,j=1,...,m.$

We can see these result in recent papers \cite{ShM2010S}, \cite{ShM2011F}. These type problems also may be easily posed in much more general contexts of $T_{\Omega } \times ...\times T_{\Omega } $ or $D\times ...\times D$ domains as follows.

We mean analytic spaces in these type domains with the quazinorms in tube or pseudoconvex domains
\[\left(\int _{D}...\left(\int _{D}\left|f\left(z_{1} ,...,z_{m} \right)\right|^{p_{1} } \delta ^{\alpha _{1} } (z_{1} )d\nu (z_{1} ) \right)^{\frac{p_2}{p_1} } ...\delta ^{\alpha _{m} } (z_{1} )d\nu (z_{m} ) \right)^{\frac{1}{p_m}} ,\]
$0<p_{j} <+\infty ,\alpha _{j} >-1,j=1,...,m,$ 

\[\left(\int _{T_{\Omega } }...\left(\int _{T_{\Omega } }\left|f\left(z_{1} ,...,z_{m} \right)\right|^{p_{1} } \Delta ^{\alpha _{1} } (z_{1} )d\tilde{\nu }(z_{1} ) \right)^{\frac{p_2}{p_1}} ...\Delta ^{\alpha _{m} } (z_{1} )d\tilde{\nu }(z_{m} ) \right)^{\frac{1}{p_m}} ,\]
$0<p_{j} <+\infty ,\alpha _{j} >-1,j=1,...,m.$

It is also rather interesting topic in function theory to consider others maps which are close and almost "diagonal" in the unit polydisk 
\[(Mf)(z_{1} ,...,z_{m} )=f(r\zeta _{1} ,...,r\zeta _{m} ),r\in (0,1),\zeta _{j} \in T,j=1,...,m;\] 
\[(Tf)(z_{1} ,...,z_{m} )=f(r_{1} \zeta ,...,r_{1} \zeta ),r_{j} \in (0,1),\zeta \in T,j=1,...,m.\] 
$z_{j} \in U=\left\{z\in {\mathbb C}:\left|z\right|<1\right\},j=1,...,m.$

We now provide below several sharp results obtained by first author in recent paper \cite{ShM2009J} (see also references there) concerning these rather simple operators in various spaces of analytic functions in the unit polydisk in $\mathbb C^{n} $.

Note similar type $M,$ $T$ operators may be considered also in more several cases in the products of unit balls and in product of tubular domains and bounded strongly pseudoconvex domains leaving then to readers also. We here omit easy details here.

Let
\[U_{*}^{n} =\left\{z\in U^{n} :\zeta \in T,r_{j} \in \left(0,1\right),j=1,...,n\right\};\] 
\[\tilde{U}^{n} =\left\{z\in U^{n} :\left|z_{j} \right|=r,r\in (0,1],j=1,...,n\right\}.\] 
\[\left\| f\right\| _{A_{\vec{\alpha }}^{p} }^{p} =\int _{U^{n} }\left|f(z)\right|^{p}  \prod _{j=1}^{n}\left(1-\left|z_{j} \right|^{2} \right)^{\alpha _{j} }  dm_{2n} (z)<+\infty ,0<p<+\infty ,\alpha _{j} >-1;\] 
\[\left\| f\right\| _{A_{\vec{\alpha }}^{p} (U_{*}^{n} )}^{p} =\int _{0}^{1}...\int _{0}^{1}\int _{T}\left|f(z)\right|^{p}    \prod _{j=1}^{n}\left(1-\left|z_{j} \right|^{2} \right)^{\alpha _{j} }  d(z_{j} )d\zeta <+\infty ,0<p<+\infty ,\alpha _{j} >-1;\] 
\[\left\| f\right\| _{A_{\alpha }^{p} (\tilde{U}^{n} )}^{p} =\int _{T^{n} }\int _{0}^{1}\left|f(z)\right|^{p}   \left(1-\left|z\right|^{2} \right)^{\alpha } d(z)d\vec{\zeta }<+\infty ,0<p<+\infty ,\alpha >-1.\] 

We have 
\[Trace\, (A_{\beta }^{p} (\tilde{U}^{n} ))=A_{\beta +n-1}^{p} (U),\]\hbox{}for $p\le 1,\beta >-1,$ and for  $p>1,\beta >\beta _{0} =\beta _{0} (n,p)$;
\[Trace\, (A_{\vec{\beta }}^{p} (U_{*}^{n} ))=A_{\left|\beta \right|+n-1}^{p} (U),\beta =(\beta _{1} ,...,\beta _{n} ),\left|\beta \right|=\sum _{k=1}^{n}\beta _{k}  ,\] 
for $p\le 1,\beta _{k} >-1,k=1,...,n,$ and for  $p>1,\beta >\beta _{0} =\beta _{0} (n,p)$.

Traces on $\left(z,...,z\right)$ of analytic spaces with the following quazinorms can also be seen in \cite{ShM2009J}. 

It will be interesting to find extensions to tube and bounded strongly pseudoconvex domains
\[\Lambda _{\alpha } (U_{*}^{n} )=\left\{f\in H(U^{n} ):\mathop{\sup }\limits_{\zeta \in T, r_{j} \in (0,1)} \left|f(r_{1} \zeta ,...,r_{n} \zeta )\right|\prod _{j=1}^{n}\left(1-r_{j} \right)^{\alpha }  <+\infty \right\},\] 
$\alpha \ge 0,j=1,...,n;$
\[\tilde{\Lambda }_{\alpha } (\tilde{U}^{n} )=\left\{f\in H(U^{n} ):\mathop{\sup }\limits_{\zeta _{j} \in T, r\in (0,1)} \left|f(r\zeta _{1} ,...,r\zeta _{n} )\right|\left(1-r\right)^{\alpha } <+\infty \right\},\] 
$\alpha \ge 0,j=1,...,n.$

Consider for example in $T_{\Omega } $ analytic spaces for $\alpha \ge 0$
\[\mathop{\sup }\limits_{\zeta _{j} \in T,r\in (0,1)} \left|f(r\zeta _{1} ,...,r\zeta _{n} )\right|\left(1-r\right)^{\alpha } <+\infty .\] 

The question to find traces of this and similar spaces in $T_{\Omega } $.

The some type question can be posed in bounded pseudoconvex domain.

Let is much more difficult issue to provide complete description of traces of various analytic function spaces in very difficult domains, namely in tubular over symmetric cones and in bounded strongly pseudoconvex domains with smooth boundary in $\mathbb C^{n} $. 

We give below concrete analytic spaces indicating that finding traces of all these analytic spaces is an open, interesting, and probably difficult issue. 

We give some quazinorms of mixed norm analytic spaces in both domains and Hardy type analytic spaces in both domains in $\mathbb C^{n} $.

Note that very interesting spaces of mixed type can be also considered in $T_{\Omega }^{m} $ or $D^{m} ,m\ge 1,$ for example the following type $\left\| \left\| \left\| f\right\| _{X_{1} } \right\| _{X_{2} } ...\right\| _{X_{m} } ,$ where $\left(X_{j} \right)=\left(A_{\alpha _{j} }^{p_{j} } \right)\left(T_{\Omega } \right)$ or $H^{p_{j} } (T_{\Omega } ),$ $j=1,...,m,$ or $\left(X_{j} \right)=\left(A_{\alpha _{j} }^{p_{j} } \right)\left(D\right)$ or $H^{p_{j} } (D),$ $j=1,...,m,$ to find estimates of traces of these analytic spaces is a very nice open problem.

For $1\le p<+\infty $ is the Hardy space  $H^{p} (T_{\Omega } )$ is the space of holomorphic functions in $T_{\Omega } $ which satisfy the estimate 

\[\left\| F\right\| _{H^{p} } =\left(\mathop{\sup }\limits_{y\in \Omega } \right)\left(\int _{{\mathbb R}^{n} }\left|F(x+iy)\right|^{p} dx \right)^{\frac{1}{p} } .\] 

In $T_{\Omega }^{m} =T_{\Omega } \times ...\times T_{\Omega } $ we can define $H^{p} (T_{\Omega }^{m} )$ in two ways with quazinorms based on quazinorms in $T_{\Omega } $
\[\left(\mathop{\sup }\limits_{y\in \Omega } \right)\left(\int _{{\mathbb R}^{n} }\left|F(x+iy)\right|^{p} dx \right)^{\frac{1}{p} } ,\] 
\[\left(\mathop{\sup }\limits_{y_{j} \in \Omega _{j} ,j=1,...,m} \right)\left(\int _{{\mathbb R}^{n} }...\int _{{\mathbb R}^{n} }\left|F(x_{1} +iy_{1} ,...,x_{m} +iy_{m} )\right|^{p} d\vec{x}  \right)^{\frac{1}{p}} \] 
or
\[\left(\mathop{\sup }\limits_{y_{m} \in \Omega } \right)\int _{{\mathbb R}^{n} }...\left(\mathop{\sup }\limits_{y_{1} \in \Omega } \right)\int _{{\mathbb R}^{n} }\left|F(x_{1} +iy_{1} ,...,x_{m} +iy_{m} )\right|^{p} dx_{1} ...dx_{m}   ,\] 
and then ask about traces of these type quazinorms for all values  $1<p<+\infty $.

This problem is open in bounded pseudoconvex domains with smooth boundary consider $\left\| \left\| \left\| f\right\| _{X_{1} } \right\| _{X_{2} } ...\right\| _{X_{m} } $function spaces quasinormed where $\left(X_{j} \right)$ is a function spaces (analytic) with one of the following quazinorms in $D$ domains
\[\sum _{\left|\alpha \right|\le k}\int _{0}^{r_{0} }\left(r^{\delta } \left(\int _{\partial D_{r} }\left|\mathcal D^{\alpha } f\right|^{p} d\sigma _{r}  \right)^{\frac{q}{p} } \frac{dr}{r} \right)^{\frac{1}{q}}   ,0<q\le +\infty ,0<p<+\infty ,k=0,1,2,...,\] 
$D_{r} =\left\{z\in D:\rho (z)>r\right\},$ $\delta >0,$ $d\sigma _{r} $ be the normalized surface measure on $\partial D_{r} $, $\mathcal D^{\alpha } $ is a fractional derivative of $f$ or 
\[\left(\mathop{\sup }\limits_{0<r<r_{0} } \right)\left(\sum _{\left|\alpha \right|\le k} \right)\left(r^{\delta } \int _{\partial D_{r} }\left|\mathcal D^{\alpha } f\right|^{p} d\sigma _{r}  \right)^{\frac{1}{p} } ,0<p<+\infty ,\delta >0,k=0,1,2,...,\] 
 $\mathop{\sup }\limits_{z\in D} \left|f(z)\right|\rho ^{\delta } (z)<+\infty ,$ $\delta \ge 0,$ $dist(z,\partial D)\, \asymp \rho (z)$

\[\left(\int _{D}\left|f(z)\right|^{p} \rho ^{\delta } (z)d\upsilon (z) \right)^{\frac{1}{p} } ,0<p<+\infty ,\delta >-1,\] 
\[\left(\int _{0}^{r_{0} }\left(\mathop{\sup \left|f\right|^{p} }\limits_{\partial D_{r} } \right)r^{\delta p-1} dr \right)^{\frac{1}{p}  } ,0<p<+\infty ,\delta >0,\] 
\[H^{p} (D_{\varepsilon } )=\left\{f\in H(D):\mathop{\sup }\limits_{\varepsilon >0} \left(\int _{\partial D_{\varepsilon } }\left|f(\zeta )\right|^{p} d\sigma (\zeta ) \right)^{\frac{1}{p}} \varepsilon ^{\alpha } <+\infty ,\alpha \ge 0,0<p<+\infty \right\},\] 
where  $\partial D_{\varepsilon } =\left\{z\in D:\rho (z)=\varepsilon \right\},$ $\varepsilon >0,$ $d\sigma $ be the normalized surface measure on $\partial D_{\varepsilon } $.

Let $B(z,r)$ be a Kobayashi ball in $D$, $z\in D,r>0.$ We can consider also Herz type quazinorms on $D$.
\[\int _{D}\left(\int _{B(z,r)}\left|f(w)\right|^{p} \delta ^{\alpha } (w)d\upsilon (w) \right)^{\frac{p}{q}  } d\upsilon (z) ,0<p<+\infty ,\alpha >-1,\] 

\[\sum _{k\ge 0}\left(\int _{B(a_{k} ,r)}\left|f(w)\right|^{q} \delta ^{\alpha } (\omega )d\upsilon (w) \right) ^{\frac{p}{q} } ,0<p,q<+\infty ,\alpha >-1,\] 
where $\left\{a_{k} \right\}$ is a known $r$-lattice in  $D$ (see \cite{ShK2015I}, \cite{ShK2016}, \cite {OrtFab1994}). All these many problems concerning   estimates of Traces of these analytic in $T_{\Omega }^{m} =T_{\Omega } \times ...\times T_{\Omega } $ or $D\times ...\times D$ are open.

We also now indicate some recent sharp results concerning traces of analytic function spaces of BMOA type and Herz spaces Bloch--type in tube, pseudoconvex domains and in the unit ball below. Some open problems in this direction will also be given by us below.

Note some results on traces of analytic function spaces of area Nevanlinna type in the polydisk can be seen in \cite{Shved1985}-\cite{DjShF1988} (and see also many references there).

We define Banach space Bloch space
\[\Lambda (r_{1} ,...,r_{k} )=\{f\in H(B^{m} ):\left\| f\right\| _{\Lambda (r_{1} ,...,r_{k} )} =\]
\[=\mathop{\sup }\limits_{z_{1} ,...,z_{k} \in B} \left|f(z_{1} ,...,z_{k} )\right|\prod _{j=1}^{n}\left(1-\left|z_{j} \right|^{2} \right)^{r_{j} } <+\infty,r_{j} >0,j=1,...,m,k\ge 1\}.\] 

\begin{theorem} Let $r_{j} >0,$ $j=1,...,m,$ $r=\sum \limits_{j=1}^{m}r_{j}$. Then $Trace\left(\Lambda (r_{1} ,...,r_{m} )\right)=\Lambda (r)$. 
\end{theorem}

Fhuter let also 
\[\Lambda _{\log } (r_{1} ,...,r_{k} )=\{f\in H(B^{m} ):\left\| f\right\| _{\Lambda (r_{1} ,...,r_{k} )} =\mathop{\sup }\limits_{z_{j} \in B} \left|f(\vec{z})\right|\prod _{j=1}^{n}\log \left(\frac{1}{1-\left|z_{j} \right|} \right)^{-\frac{1}{r_j}} <+\infty \},\]
$\displaystyle \sum \limits_{j=1}^{m}\frac{1}{r_{j} } =1,$ $r_{j} >0.$

Then we have that 
\[Trace\, \left(\Lambda _{\log } (r_{1} ,...,r_{m} )\right)=\left\{f\in H(B):\mathop{\sup }\limits_{z\in B} \left|f(z)\right|\left(\log \frac{1}{1-\left|z\right|} \right)^{-1} \left(1-\left|z\right|\right)<+\infty \right\}.\] 

It will be interesting to find complete analogues of these results in tube, pseudoconvex domains.

We provide traces of BMOA in the ball, in pseudoconvex and tube domains similar result can be seen \cite{ShM2017C}, \cite{ShM2022R} but $E(0)=B,$ $a\in B\backslash \{ 0\} ,$ $E(a)=\left\{z\in B:\left|1-\left\langle z,\frac{a}{\left|a\right|} \right\rangle \right|<1-\left|a\right|\right\}.$

Let $k\ge 1,$ $k\in {\mathbb Z},$ $\mu $ is a positive Borel measure on $\mathop{B\times ...\times B} ,$ $r_{j} >0,$ $j=1,...,k,$ $\mu $ is a bounded $\left(r_{1} ,...,r_{k} \right)$ Carleson measure if
\[\sup \frac{\mu \left(E(a_{1} )\times ...E(a_{k} )\right)}{\prod \limits_{j=1}^{k}\left(1-\left|a_{j} \right|\right)^{r_{j} }  } <+\infty ,a_{j} \in B,j=1,...,k.\] 

Let $0<p<+\infty ,s_{1} ,...,s_{k} >-1,r_{1} ,...,r_{k} >0$, $M_{r_{1} ,...,r_{k} }^{} \left(B^{k} ;d\nu _{s_{1} } ,...,d\nu _{s_{m} } \right)$ is a space of measurable functions on $B^{k} $ for $\mu _{f} =\left|f\right|^{p} \prod\limits _{j=1}^{k}d\nu _{s_{j} }  $ is a bounded $\left(r_{1} ,...,r_{m} \right)$ Carleson measure with
\[\sup \frac{\mu _{f} (E(a_{1} ),...,E(a_{m} ))}{\prod \limits_{j=1}^{k}\left(1-\left|a_{j} \right|\right)^{r_{j} }  } ,a_{j} \in B,j=1,...,k.\] 

\begin{theorem} Let $t=(m-1)(n+1)+\sum \limits_{j=1}^{m}s_{j}  $. 

Then $Trace\left(H(B^{m} )\cap M_{r_{1} ,...,r_{m} }^{p} (B^{m} ;d\nu _{s_{1} } ,...,d\nu _{s_{m} } )\right)\, =H(B)\cap M_{r}^{p} (B;d\nu _{t} ),$ $0<p<+\infty ,$ $s_{j} >-1,$ $r_{j} >0,$ $j=1,...,m,$ $r=\sum \limits_{j=1}^{m}r_{j}  ,$ $r_{j} <n+1+s_{j} ,$ $j=1,...,m.$
\end{theorem}

Some descriptions (or estimates) of traces of analytic spaces $\left(M_{\alpha }^{p} \right),$ $\left(K_{\alpha ,\beta }^{p,q} \right), $ $\left(D_{\alpha ,\beta }^{p} \right)$ and related spaces 
\[\left(M_{\alpha }^{p} \right)=\left\{f\in H(B^{m} ):\int _{S}\int _{\Gamma _{t} (\zeta )}...\int _{\Gamma _{t} (\zeta )}\left|f(z)\right|^{p} d\nu _{\alpha } (z)d\sigma (\zeta )<+\infty    \right\},\]
$0<p<+\infty ,$ $\alpha >-1;$
\[\left(K_{\alpha ,\beta }^{p,q} \right)=\left\{f\in H(B^{m} ):\int _{B}...\int _{B}\left(\int _{D(z_{1} ,r)}...\int _{D(z_{m} ,r)}\left|f(z)\right|^{p} d\nu _{\alpha } (z)  \right)^{\frac{q}{p} } d\nu _{\beta } (z)<+\infty   \right\},\]
$0<p,q<+\infty ,$ $\alpha >-1,$ $\beta >-1;$\\
$\left(D_{\alpha ,\beta }^{p} \right)=\Bigg\{f\in H(B^{m} ):$
\[\int _{0}^{1}(1-r)^{\beta }  \left(\int _{\left|z_{1} \right|<r}...\int _{\left|z_{m} \right|<r}\left|f(z)\right|^{p} \prod _{j=1}^{m}\left(1-\left|z_{j} \right|\right)^{\alpha _{j} }  d\nu (z_{j} )  \right)dr<+\infty \Bigg\},\] 
$0<p<+\infty ,$ $\beta >-1,$ $\alpha _{j} >-1,$ $j=1,...,m.$

Many related and (or) very close to these trace theorems interesting  results and interesting  problems were considered also by various authors in recent decades in various analytic spaces in various domains in $\mathbb C^n,$ we mention, for example, papers \cite{HenPol1984}-\cite{Am1979}.

We will not discuss these problems or these  type (trace type) problems in this expository paper referring the authors to these mentioned interesting papers.

In practically all proofs of various trace theorems, various Bergman projection type theorems are very important, and we refer the interested readers to \cite{ShTom2023}-\cite{ShMakh2024} for various new interesting results in this direction.

Our sharp results on traces in harmonic and analytic function spaces of several variables which we provided in this paper were recently extended to various spaces of pluriharmonic functions of several variables of Bergman type. We refer the reader to \cite{ShPov2013} for these results.

We finally wish to refer the reader to \cite{ShTom2021R}-\cite{ShK2014J} for many new recent interesting results of the first author and his various coauthors on embeddings in analytic Bergman type spaces  and on Bergman projections in tubular domains over cones and in bounded strongly pseudoconvex domains with smooth boundary which can be applied in the future in various ways or used to get proofs of new trace theorems in various analytic function spaces in such types complicated domains in $\mathbb C^n$.

We refer the reader to \cite{Hen1972}, \cite{Cum1983} and many references there for some old interesting results closely related to our trace theorems.

Probably, various embeddings for analytic $Q_p$ spaces from \cite{Xi2001} can be used in proofs of trace theorems for such type spaces in the unit disk and polydisk.

We refer the reader to \cite{Ren2004} and many references there for some new results on sharp traces in simple polydisk in mixed norm spaces.

Some embeddings from \cite{Fr1972} can be used to prove traces theorems in various  Hardy type spaces in the polydisk.

Many embeddings from \cite{{Zhu2005}} can be generalized to bounded strongly pseudoconvex domains with smooth boundary and then used to get trace theorems in various analytic spaces in such domains, we pose this as a problem for interested readers.

\vskip 1 cm \footnotesize
\begin{flushleft}
Romi F.~Shamoyan, Nataliya M. Makhina \\ 
 Department of Mathematics and Physics\\ 
 Bryansk State Academician I.G. Petrovski University\\ 
 14 Bezhitskaya St,\\ 
241036 Bryansk, Russia \\ 
E-mails: rsham@mail.com, makhina32@gmail.com 
\end{flushleft}

\vskip0.5cm
\begin{flushright}
Received: ??.??.2025
\end{flushright}

 \end{document}